\documentclass[11pt, reqno]{amsart}

\usepackage{amssymb}
\usepackage{amsmath}
\usepackage{enumerate}

\def\stackunder#1#2{\mathrel{\mathop{#2}\limits_{#1}}}

\makeatletter\@addtoreset{equation}{section}\makeatother

\allowdisplaybreaks

\begin{document}

\bibliographystyle{alpha}

\newtheorem{theorem}{Theorem}[section]
\newtheorem{proposition}[theorem]{Proposition}
\newtheorem{definition}[theorem]{Definition}
\newtheorem{corollary}[theorem]{Corollary}
\newtheorem{lemma}[theorem]{Lemma}
\newtheorem{remark}[theorem]{Remark}
\newtheorem{example}[theorem]{Example}
\newtheorem{conjecture}[theorem]{Conjecture}

\title[N/V-limit for stochastic dynamics]{N/V-limit for
stochastic dynamics in continuous particle systems}
\author{Martin Grothaus, Yuri G.~Kondratiev, Michael R\"ockner}

\address{Martin Grothaus, Mathematics Department,
University of Kaiserslautern, \newline P.O.Box 3049, 67653
Kaiserslautern, Germany; BiBoS, Bielefeld University, 33615
Bielefeld, Germany and SFB 611, IAM, University of Bonn, 53115
Bonn, Germany.
\newline {\rm \texttt{Email: grothaus@mathematik.uni-kl.de},
\texttt{URL: http://www.mathematik.uni-kl.de/$\sim$grothaus/ }}
\newline
Yuri G.~Kondratiev, BiBoS and Mathematics Department,
Bielefeld University, 33615 Bielefeld, Germany and
Inst.~Math., NASU, 252601 Kiev, Ukraine.
\newline
{\rm \texttt{Email: kondrat@mathematik.uni-bielefeld.de}}.
\newline
Michael R\"ockner, BiBoS and Mathematics Department,
Bielefeld University, 33615 Bielefeld, Germany. \newline
{\rm \texttt{Email: roeckner@mathematik.uni-bielefeld.de}
}}

\date{\today}

\thanks{
\\
We thank Hans-Otto Georgii, Oleksandr Kutoviy, Robert Minlos 
and Herbert Spohn for discussions and helpful comments. 
Financial support through CCM, University of
Madeira, is gratefully acknowledged.}


\begin{abstract} We provide an $N/V$-limit for the
infinite particle, infinite volume stochastic dynamics associated
with Gibbs states in continuous particle systems on ${\mathbb
R}^d, d \ge 1$. Starting point is an $N$-particle stochastic
dynamic with singular interaction and reflecting boundary
condition in a subset $\Lambda \subset {\mathbb R}^d$ with finite
volume (Lebesgue measure) $V = |\Lambda| < \infty$. The aim is to
approximate the infinite particle, infinite volume stochastic
dynamic by the above $N$-particle dynamic in $\Lambda$ as $N \to
\infty$ and $V \to \infty$ such that $N/V \to \rho$, where $\rho$
is the particle density. First we derive an improved Ruelle bound
for the canonical correlation functions under an appropriate
relation between $N$ and $V$. Then tightness is shown by using the
Lyons--Zheng decomposition. The equilibrium measures of the
accumulation points are identified as infinite volume canonical
Gibbs measures by an integration by parts formula and the
accumulation points themselves are identified as infinite
particle, infinite volume stochastic dynamics via the associated
martingale problem. Assuming a property closely related to 
Markov uniqueness and weaker than essential self-adjointness,
via Mosco convergence techniques we can identify 
the accumulation points as Markov processes
and show uniqueness. I.e., all accumulation corresponding 
to one invariant canonical Gibbs measure coincide.  
The proofs work for general repulsive
interaction potentials $\phi$ of Ruelle type and all temperatures,
densities, and dimensions $d \ge 1$, respectively.
$\phi$ may have a nontrivial negative
part and infinite range as e.g.~the Lennard--Jones potential.
Additionally, our result provides as a by-product an approximation
of grand canonical Gibbs measures by finite volume canonical Gibbs
measures with empty boundary condition.
\end{abstract}

\maketitle

\section{Introduction}
The infinite particle, infinite volume stochastic dynamics $({\bf
X}(t))_{t \ge 0}$ in continuous particle systems is an infinite
dimensional diffusion process having a Gibbs measure $\mu$,
e.g.~of the type studied by Ruelle in \cite{Rue69}, as an
invariant measure. Physically, it describes the stochastic
dynamics of infinite Brownian particles in ${\mathbb R}^d, \, d
\ge 1,$ which are interacting via the gradient of a pair-potential
$\phi$. Since each particle can move through each position in
space, the system is called continuous and is used for modelling
gas and fluids. For realistic models which can be described by
these stochastic dynamics, e.g.~suspensions, we refer to
\cite{Sp86}.

The infinite particle, infinite volume stochastic dynamics takes
values in the configuration space
\begin{eqnarray*}
\Gamma :=\big\{\,\gamma\subset {\mathbb R}^d\mid \#(\gamma \cap
\Lambda)<\infty\text{ for each compact }\Lambda\subset {\mathbb
R}^d\,\big\},
\end{eqnarray*}
and informally solves the following infinite system of stochastic
differential equations:
\begin{align}\label{infinitvoldyn}
dx(t) & = -\beta \!\! \sum_{\stackunder{y(t) \neq x(t)} {y(t) \in
{\bf X}(t)}} \nabla \phi(x(t) - y(t)) \,dt + \sqrt{2}\,dB^x(t), &
 \nonumber \\
{\bf P} \circ {\bf X}(0)^{-1} & = \mu,
\end{align}
where $x(t) \in {\bf X}(t) \in \Gamma$, $(B^{x})_{x \in \gamma}$,
$\gamma \in \Gamma$, is a sequence of independent Brownian
motions and $\mu$ is the
invariant measure. The study of such diffusions has been initiated
by R.~Lang \cite{La77} (see also \cite{Sh79}), who considered the
case $\phi \in C_0^3({\mathbb R}^d)$ using finite dimensional
approximations and stochastic differential equations. More
singular $\phi$, which are of particular interest in Physics, as
e.g.~the Lennard--Jones potential, have been treated by H.~Osada,
\cite{Os96}, and M.~Yoshida, \cite{Y96} (see also \cite{Ta97},
\cite{FRT00} for the hard core case). Osada and Yoshida were the
first to use Dirichlet forms for the construction of such
processes. However, they could not write down the corresponding
generators or martingale problems explicitly, hence could not
prove that their processes actually solve (\ref{infinitvoldyn})
weakly. This, however, was proved in \cite{AKR98b} by showing an
integration by parts formula for the respective Gibbs measures. In
\cite{AKR98b}, also Dirichlet forms were used and all
constructions were designed to work particularly for singular
potentials of the above mentioned type. Additionally, an explicit
expression for the corresponding generator and martingale problem
was provided, which shows that the process in \cite{AKR98b} indeed
solves \eqref{infinitvoldyn} in the weak sense.

In this paper, by an approximation through $N$-particle stochastic
dynamics in subsets $\Lambda \subset {\mathbb R}^d$ with finite
volume (Lebesgue measure) $V = |\Lambda| < \infty$, we construct
weak solutions to \eqref{infinitvoldyn}. The approximation is done
in terms of the $N/V$ limit, i.e., $N \to \infty$ and $V \to
\infty$ such that $N/V \to \rho$, where $\rho$ is the particle
density.

The $N$-particle stochastic dynamics in $\Lambda$, $({\bf
X}(t))_{t \ge 0}$, takes values in the space of $N$-point
configurations in $\Lambda$:
\begin{eqnarray*}
\Gamma^{\scriptscriptstyle{(N)}}_{\scriptscriptstyle{\Lambda}} :=
\{ \gamma \subset \Lambda \,|\, \# (\gamma) = N \} \subset \Gamma.
\end{eqnarray*}
It solves weakly the following $N$-system of stochastic
differential equations before hitting $\partial
(\Gamma_{\scriptscriptstyle{\Lambda}}^{\scriptscriptstyle{(N)}})$:
\begin{align}
dx(t) & = - \beta \!\! \sum_{\stackunder{y(t) \, \neq \, x(t)}
{y(t) \in \,\, {\bf X}(t)}}
\nabla \phi(x(t) - y(t)) \,dt + \sqrt{2} \,dB^{x_0}(t), \nonumber \\
& \mbox{with reflecting boundary condition}, \label{npsdintro}
\end{align}
for sufficiently many initial conditions $\gamma_0 \in
\Gamma_{\scriptscriptstyle{\Lambda}}^{\scriptscriptstyle{(N)}}$.
Here $x(t) \in {\bf X}(t) \in
\Gamma_{\scriptscriptstyle{\Lambda}}^{\scriptscriptstyle{(N)}}$
and $(B^{x_0})_{x_0 \in \gamma_0}$ are $N$ independent Brownian
motions starting in $x_0$. A weak solution to \eqref{npsdintro}
has been constructed in \cite{FaGr04a}, see Theorem \ref{exnpsd}.
There the authors have used the Dirichlet form approach and their
construction works for all dimensions and very general interaction
potentials $\phi$. Essentially, the interaction potential $\phi$
only has to have a singularity at the origin (repulsion) (RP), to
be bounded from below (BB) and weakly differentiable (D), see
below for a precise definitions.

Note that we are only considering configurations with at most one
particle in one position, which is a reasonable assumption for
modelling gas and fluids. In such a setting for dimension $d=1$
this is the first existence result for a solution to
\eqref{npsdintro}. The essential assumption for this result is the
condition (RP) (repulsion of close particles), which is natural
from the physical point of view.


Our approach is different from the finite dimensional
approximation provided by Lang \cite{La77}. There for a fixed
subset $\Lambda \subset {\mathbb R}^d$ with finite volume, the
finite particle, finite volume dynamics consists of finitely but
arbitrarily (for different initial conditions) many interacting 
particles inside the volume and additionally they are interacting with particles from
the complement of $\Lambda$. That construction is rather in a
grand canonical setting whereas ours is in a canonical one.
Thus, we expect the finite particle, finite volume dynamics used in \cite{La77}
for singular interaction potentials with non-trivial negative part 
not to have such nice properties as our
$N$-particle stochastic dynamics in $\Lambda$. E.g., 
for determining a spectral gap of their generators, it is much nicer
to have a fixed number $N$ of particles in a given volume $\Lambda$
not interacting with particles in the complement of $\Lambda$, than
finite but arbitrarily many particles inside $\Lambda$ interacting
with in general infinitely many particles in the complement of $\Lambda$.  

Our plan for future work is to use our approximation by nice
processes to get better knowledge about the infinite volume,
infinite particle dynamics. For example we would like to: explore
in more detail the structure of the spectrum of its generator and
study the problem of essential self-adjointness; construct
non-equilibrium infinite particle, infinite volume stochastic
dynamics; tackle the Boltzmann--Gibbs principle, see
e.g.~\cite{Sp86} and \cite{GKLR00}; use our approximation technique
to construct solutions to other equations as e.g.~the Langevin equations.

The present paper is organized in the following way: In Section
\ref{s1} we define a metric on the configuration space $\Gamma$
which is appropriate for our problem. This metric is from the
class of metrics on $\Gamma$ developed in \cite{KoKu04} and
induces the vague topology. Essential for our considerations is
that these metrics $d$ make $(\Gamma, d)$ a Polish space and that
relative compact sets w.r.t.~the vague topology can be described
explicitly (cf.~\cite{KoKu04}).

The concept of canonical Gibbs measures, our assumptions on the
interaction potential and a precise definition of the $N/V$-limit
are presented in Section \ref{s2}. Furthermore, in Theorem
\ref{ruelle} we prove the first major result of this paper. There
we establish a bound for canonical correlation functions analogous
to the Ruelle bound for grand canonical correlation functions, see
\eqref{Ruellebo}. In the proof we combine ideas of Ruelle's proof
\cite{Rue70} for deriving the Ruelle bound in the grand canonical
case with estimates obtained in \cite{DoMi67}. A major difference
in comparison with the grand canonical case is that in the
canonical case a right balance between the particle number $N$ and
the volume $V$ is necessary. Furthermore, we derive an improved
Ruelle bound for canonical correlation functions, see
\eqref{expRuellebo}. This bound enables us to take into account
potentials with singularities at the origin, see condition (D)
below.

In Section \ref{s3} we briefly summarize the construction
of the $N$-particle stochastic dynamics in $\Lambda$
weakly solving \eqref{npsdintro} provided in \cite{FaGr04a}.

The $N/V$-limit of $N$-particle, finite volume stochastic dynamics
is then derived in Section \ref{s4}. First, in Theorem \ref{tight}
we prove tightness of the sequence of laws $({\bf
P}^{\scriptscriptstyle{(N)}})_{N \in {\mathbb N}}$ of the
equilibrium $N$-particle, finite volume stochastic dynamics in the
$N/V$-limit. Equilibrium stochastic dynamics means that the
stochastic dynamics starts with an initial distribution given by
the corresponding invariant, finite volume canonical Gibbs measure
${\mu}^{\scriptscriptstyle{(N)}}$. The proof is split into two
lemmas. Lemma \ref{tight0} gives tightness of the corresponding
one-dimensional distributions (invariant, finite volume canonical
Gibbs measures) $({\mu}^{\scriptscriptstyle{(N)}})_{N \in {\mathbb
N}}$ and essentially depends on the improved Ruelle bound
\eqref{expRuellebo} and the description of compact sets provided
in \cite{KoKu04}. In Lemma \ref{tight1} we prove
Kolmogorov--Chentsov type estimates for the increments. In the
proof we use the well-known Lyons--Zheng decomposition,
\cite{LZ88}, \cite{LZ94}, of the $N$-particle, finite volume
stochastic dynamics and the Burkholder--Davies--Gundy inequalities
in order to establish the required estimate of the increments. For
this it is important to have sufficiently many functions in the
domain of the corresponding Dirichlet form, which is in fact
implied by the reflecting boundary condition we impose on the
$N$-particle, finite volume stochastic dynamics. Again, also the
improved Ruelle bound is of essential importance.

Then in Theorem \ref{thidibp} we prove an integration by parts
formula for the accumulation points $\mu$ of
$({\mu}^{\scriptscriptstyle{(N)}})_{N \in {\mathbb N}}$. Together
with a characterization theorem provided in \cite{AKR98b} this
implies that these $\mu$ are infinite volume canonical Gibbs
measures.

After that, in Theorem \ref{martingaleproblem} we identify the
accumulation points ${\bf P}$ of $({\bf
P}^{\scriptscriptstyle{(N)}})_{N \in {\mathbb N}}$ as solutions of
\eqref{infinitvoldyn} in the sense of the associated martingale
problem. See also Remark \ref{solution}. In the proof we are using
that the $N$-particle, finite volume stochastic dynamics solves
the martingale problem corresponding to \eqref{npsdintro}.

From Theorem \ref{martingaleproblem} we can not conclude that
the accumulation points ${\bf P}$ of $({\bf
P}^{\scriptscriptstyle{(N)}})_{N \in {\mathbb N}}$ are laws 
of Markov processes. However, assuming a property closely related to 
Markov uniqueness and weaker than essential self-adjointness,
in Theorem \ref{thmmoscoconv} we can show Mosco convergence, 
\cite{Mos94}, \cite{KuSh03}, of 
the quadratic forms (Dirichlet forms) corresponding 
to convergent subsequences. This implies strong convergence 
of the associated semi-groups. 
This convergence, in turn, enables us
to identify the accumulation points ${\bf P}$ as laws of Markov processes 
and show uniqueness. I.e., all accumulation points ${\bf P}$
corresponding to one invariant canonical Gibbs measure coincide, see 
Theorem \ref{uniqueness}.   

Finally, in Section \ref{s5} we apply our results to the problem
of equivalence of ensembles. More precisely, as a by-product of
the results described above we obtain an approximation of grand
canonical Gibbs measures by finite volume canonical Gibbs measures
with empty boundary condition, see Theorem \ref{coroequildht}.

The progress achieved in this paper may be summarized by the
following list of main results:
\begin{itemize}
\item Derivation of an improved Ruelle bound for canonical correlation
functions, see Theorem \ref{ruelle}.

\item Tightness of the sequence of laws $({\bf
P}^{\scriptscriptstyle{(N)}})_{N \in {\mathbb N}}$ of
equilibrium $N$-particle, finite volume stochastic dynamics in the
$N/V$-limit, see Theorem \ref{tight}.

\item Identification of the accumulation points $\mu$ of the
sequence of finite volume canonical Gibbs measures
$({\mu}^{\scriptscriptstyle{(N)}})_{N \in {\mathbb N}}$ as
infinite volume canonical Gibbs measures via an integration by
parts formula, see Theorem \ref{thidibp}.

\item Identification of the accumulation points ${\bf P}$
of the sequence of laws $({\bf P}^{\scriptscriptstyle{(N)}})_{N \in {\mathbb N}}$ 
of equilibrium $N$-particle, finite volume stochastic dynamics in the
$N/V$-limit as solutions of \eqref{infinitvoldyn} in the sense of
the associated martingale problem, see Theorem
\ref{martingaleproblem}. This is the first construction of a
solution to \eqref{infinitvoldyn} for $d=1$ with state space
$\Gamma$ (at most one particle in one position).
\end{itemize}

Furthermore, when assuming a property closely related to 
Markov uniqueness and weaker than essential self-adjointness:

\begin{itemize}
\item Identification of the accumulation points ${\bf P}$
of the sequence of laws $({\bf
P}^{\scriptscriptstyle{(N)}})_{N \in {\mathbb N}}$ of
equilibrium $N$-particle, finite volume stochastic dynamics in the
$N/V$-limit as Markov processes and showing uniqueness, see 
Theorem \ref{uniqueness}.
\end{itemize}

At the moment we are working on the assumed property and expect
to show it soon.
 
All above results apply to all dimensions $d \ge 1$,
temperatures and densities and to physically relevant repulsive
(RP) interaction potentials $\phi$. Additional assumptions are
only a mild temperedness (T) condition (fast enough decay in the
long range), that the potential is bounded from below (BB) and a
mild differentiability (D) condition. Hence, singularities at the
origin, non-trivial negative part, and infinite range are allowed.

Hypotheses on the potential are weakened not for the sake of
generality, but in order to cover the physically relevant
potentials (as e.g.~ Lennard--Jones potential).

\section{A Polish metric for the configuration space}\label{s1}

The configuration space $\Gamma$ over ${\mathbb R}^d$,
$d\in{\mathbb N}$, is defined as the set of all subsets of
${\mathbb R}^d$ which are locally finite:
\begin{eqnarray*}
\Gamma :=\big\{\,\gamma\subset {\mathbb R}^d\mid
\#(\gamma_\Lambda)<\infty\text{ for each compact }\Lambda\subset
{\mathbb R}^d\,\big\},
\end{eqnarray*}
where $\#$ denotes the number of elements of a set and
$\gamma_\Lambda:= \gamma\cap\Lambda$. One can identify
$\gamma\in\Gamma$ with the positive Radon measure
$\sum_{x\in\gamma}\varepsilon_x\in{\mathcal M}({\mathbb R}^d)$,
where  $\varepsilon_x$ is the Dirac measure at $x$,
$\sum_{x\in\varnothing}\varepsilon_x :=$ zero measure, and
${\mathcal M}({\mathbb R}^d)$ stands for the set of all positive
Radon  measures on the Borel $\sigma$-algebra ${\mathcal
B}({\mathbb R}^d)$. A metric on ${\mathcal M}({\mathbb R}^d)$ is
given by
\begin{eqnarray*}
d_{\mathcal M}(\nu, \mu) := \sum_{k = 1}^\infty 2^{-k} p_k \left( 1-
\exp \left(-|\langle f_k, \nu-\mu \rangle| \right) \right),
\qquad \nu, \mu \in {\mathcal M}({\mathbb R}^d),
\end{eqnarray*}
where $\{f_k \,|\, k \in {\mathbb N} \} \subset C^1_c({\mathbb
R}^d)$ (space of continuously differentiable functions on
${\mathbb R}^d$ with compact support) is a measure determining
class, $(p_k)_{k \in {\mathbb N}}$ a sequence of strictly positive
weights bounded by $1$, and
\begin{eqnarray*}
\langle f, \nu \rangle := \int_{{\mathbb R}^d} f \,d\nu, \quad f
\in C_c({\mathbb R}^d), \,\, \nu \in {\mathcal M}({\mathbb R}^d).
\end{eqnarray*}
$\{f_k \,|\, k \in {\mathbb N} \}$ can be chosen so that
$d_{\mathcal M}$ induces the vague topology on ${\mathcal
M}({\mathbb R}^d)$. This metrization is  separable and complete,
see \cite[A~7.7]{Ka75} for the case $\{ f_k \,|\, k \in {\mathbb
N}\} \subset C_c({\mathbb R}^d)$ and $p_k = 1$ for all $k \in
{\mathbb N}$.

In ${\mathcal M}({\mathbb R}^d)$ we consider the subset ${\mathcal
R}({\mathbb R}^d)$ consisting of all ${\mathbb Z}_+\cup
\{\infty\}$-valued Radon measures. Since ${\mathcal R}({\mathbb
R}^d)$ is a closed subset of ${\mathcal M}({\mathbb R}^d)$
w.r.t.~the vague convergence, see \cite[A~7.4]{Ka75}, also
$({\mathcal R}({\mathbb R}^d), d_{\mathcal M})$ is a Polish space.

Now our aim is to find a metric on $\Gamma$ which is Polish. Let
$\Phi: (0, \infty) \to [0, \infty)$ be a continuous decreasing
function such that $\lim_{t \to 0} \Phi(t) = \infty$; and let $h:
{\mathbb R}^d \to (0,1]$ be a function in $L^1({\mathbb R}^d) \cap
C^1({\mathbb R}^d)$. Furthermore, let ${\bf I} = \{I_k \,|\, k \in
{\mathbb N}\}$ be a collection of functions from $C^1_c({\mathbb
R}^d)$ such that $I_k: {\mathbb R}^d \to [0,1]$, supp$I_k \subset
B_k(0)$, and $I_{k+1}(x) = 1$ for all $x \in B_k(0)$, here
$B_k(0)$ denotes the closed ball with radius $k$ centered at the
origin. Define
\begin{eqnarray*}
S^{\Phi,f}(\gamma) := \sum_{\{x,y\} \subset \gamma}
\exp(\Phi(|x-y|)) f(x)f(y),
\end{eqnarray*}
where $f : {\mathbb R}^d \to [0, \infty)$ is a continuously
differentiable function.
For any $k \in {\mathbb N}$ set $h_k := h I_k$.
Then for $\gamma, \eta \in \Gamma$ we define the metric
\begin{eqnarray}\label{metric}
d_{\Phi,h} (\gamma, \eta) := d_{{\mathcal M}}(\gamma, \eta) +
\sum_{k=1}^\infty 2^{-k} q_k \frac{\left| S^{\Phi,h_k}(\gamma)-
S^{\Phi,h_k}(\eta)\right|}{1 + \left| S^{\Phi,h_k}(\gamma)-
S^{\Phi,h_k}(\eta)\right|},
\end{eqnarray}
where $(q_k)_{k \in {\mathbb N}}$ is a sequence of strictly positive
weights bounded by $1$. The following has been proved in \cite[Theo.~3.5, Prop.~3.1]{KoKu04} for
$q_k=1$, $k \in {\mathbb N}$. Easily, the statement generalizes to the present
situation.
\begin{proposition}\label{propmetric}
$(\Gamma, d_{\Phi,h})$ is a complete and separable metric space.
Moreover, the topology on $\Gamma$ generated by the metric
$d_{\Phi,h}$ is equivalent to the vague topology on $\Gamma$ and
the sets
\begin{eqnarray*}
\{\gamma \in \Gamma \,|\, S^{\Phi,h}(\gamma) \le R \}, \qquad R <
\infty,
\end{eqnarray*}
are relative compact subsets w.r.t.~the vague topology.
\end{proposition}

\section{Canonical Gibbs measures and an improved Ruelle bound}\label{s2}

Let $\Lambda\subset {\mathbb R}^d$. We denote
$\Gamma_{\scriptscriptstyle{\Lambda}} :=
\{\gamma\in\Gamma\mid\gamma\subset\Lambda\}$. For any $N \in
{\mathbb N}$ and bounded Borel measurable $\Lambda \subset
{\mathbb R}^d$ we define the space of $N$-point configurations in
$\Lambda$ by
\begin{eqnarray*}
\Gamma^{\scriptscriptstyle{(N)}}_{\scriptscriptstyle{\Lambda}} :=
\{ \gamma \subset \Lambda \,|\, \# (\gamma) = N \} \subset
\Gamma_{\scriptscriptstyle{\Lambda}}.
\end{eqnarray*}
To define more structure on
$\Gamma^{\scriptscriptstyle{(N)}}_{\scriptscriptstyle{\Lambda}}$
we use the following natural mapping
\begin{align*}
\mbox{sym}^{\scriptscriptstyle{(N)}}: \widetilde{\Lambda^N} &\to
\Gamma^{\scriptscriptstyle{(N)}}_{\scriptscriptstyle{\Lambda}} \\
\mbox{sym}^{\scriptscriptstyle{(N)}}((x_1, \ldots, x_N)) &:=
\{x_1, \ldots, x_N \},
\end{align*}
where
\begin{eqnarray*}
\widetilde{\Lambda^N} := \{(x_1, \dots, x_N) \in \Lambda^N \,|\,
x_k \neq x_j \quad \mbox{if} \quad k \neq j\}.
\end{eqnarray*}
These mappings generate a topology and corresponding Borel
$\sigma$-algebra on
$\Gamma^{\scriptscriptstyle{(N)}}_{\scriptscriptstyle{\Lambda}}$.
Obviously, this $\sigma$-algebra coincides with the Borel
$\sigma$-algebra inherited from $\Gamma$ equipped with its vague
topology. We denote by $dx_{\scriptscriptstyle{\Lambda}}$ the
Lebesgue measure on $\Lambda$. Then the product measure
$dx_{\scriptscriptstyle{\Lambda}}^{\otimes N}$ can be considered
on $\widetilde{\Lambda^N}$. Let
$dx_{\scriptscriptstyle{\Lambda}}^{\scriptscriptstyle{(N)}} :=
dx_{\scriptscriptstyle{\Lambda}}^{\otimes N} \circ
(\mbox{sym}^{\scriptscriptstyle{(N)}})^{-1}$ be the corresponding
measure on
$\Gamma^{\scriptscriptstyle{(N)}}_{\scriptscriptstyle{\Lambda}}$.


A pair potential (without hard core) is a Borel measurable
function $\phi\colon {\mathbb R}^d \to {\mathbb R} \cup \infty$
such that $\phi(-x)=\phi(x)\in{\mathbb R}$ for all $x\in{\mathbb
R}^d\setminus\{0\}$. For bounded Borel measurable $\Lambda \subset
{\mathbb R}^d$ the potential energy $E_{\phi}: \Gamma_{\scriptscriptstyle{\Lambda}}
\to {\mathbb R}$ in $\Lambda$ with empty boundary condition is
defined by
\begin{eqnarray*}
E_{\phi}(\gamma) := \sum\limits_{ \{x,y\} \subset \gamma}
\phi(x-y), \qquad \gamma \in \Gamma_{\scriptscriptstyle{\Lambda}},
\end{eqnarray*}
where the sum over the empty set is defined to be zero. The
interaction energy between two configurations $\gamma$ and $\eta$
from $\Gamma_{\scriptscriptstyle{\Lambda}}$ is defined by
\begin{eqnarray*}
W_\phi(\gamma, \eta) := \sum_{x\in\gamma,\, y\in\eta }\phi(x-y).
\end{eqnarray*}
Note that
\begin{eqnarray*}
E_{\phi}(\gamma \cup \eta) =  E_{\phi}(\gamma) + W_\phi(\gamma,
\eta) + E_{\phi}(\eta) , \qquad \gamma, \eta \in \Gamma_{\scriptscriptstyle{\Lambda}}.
\end{eqnarray*}

Now we fix our assumptions on $\phi$:
\begin{description}
\item[(RP)] ({\it Repulsion}) There exists a decreasing continuous
function $\Phi: (0, \infty) \to [0, \infty)$ with $\lim_{t \to 0}
\Phi(t) t^d = \infty$ and $R_1 > 0$ such that
\begin{eqnarray*}
\phi(x) \ge \Phi(|x|) \quad \mbox{for} \quad |x| \le R_1.
\end{eqnarray*}
Furthermore, the potential $\phi$ is bounded from above on $\{x
\in {\mathbb R}^d \,|\, r \le |x|\le R_1\}$ for all $r > 0$.

\item[(T)] ({\it Temperedness}) The exists $A, R_2 < \infty$ and
$\lambda > d$ such that
\begin{eqnarray*}
|\phi(x)| \le A|x|^{-\lambda} \quad \mbox{for} \quad |x| \ge R_2.
\end{eqnarray*}

\item[(BB)] ({\it Bounded below})
 There exist $B\ge0$ such that
\begin{eqnarray*}
\phi(x) \ge - B, \quad \mbox{for all} \quad x \in {\mathbb R}^d.
\end{eqnarray*}
\end{description}

For every $r=(r^1,\dots,r^d)\in{\mathbb Z}^d$, we define a cube
\begin{eqnarray*}
Q(r):=\left\{\, x\in{\mathbb R}^d \, \Bigg| \, r^i-\frac{1}{2}\le
x^i<r^i+\frac{1}{2} \,\right\}.
\end{eqnarray*}
These cubes form a partition of ${\mathbb R}^d$. For any
$\gamma\in\Gamma$, we set $\gamma_r:=\gamma_{Q(r)}$, $r\in{\mathbb
Z}^d$.

(RP), (T) and (BB) imply that $\phi$ is superstabile (SS) and
lower regular (LR), see \cite[Prop.~1.4]{Rue70}. That is:

\begin{description} \item[(SS)] ({\it Superstability})
There exist $D>0$, $K\ge0$ such that, if $\gamma\in\Gamma_
{\Lambda}$, where $\Lambda$ is a finite union of the cubes $Q(r)$,
then
\begin{eqnarray*}
\sum_{\{x,y\}\subset\gamma}\phi(x-y)\ge\sum_{r\in{\mathbb
Z}^d}\big(D\#(\gamma_r)^2-K\#(\gamma_r)\big).
\end{eqnarray*}

\item[(LR)] ({\it Lower regularity}) There exists a decreasing
positive function $\Psi: {\mathbb N} \to [0, \infty)$ such that
\begin{eqnarray*}
\sum_{r\in{\mathbb Z}^d}\Psi(|r|_{\max})<\infty
\end{eqnarray*}
and for any disjoint $\Lambda', \Lambda''$ which are finite unions
of the cubes $Q(r)$, we have for $\gamma'\in\Gamma_{\Lambda'}$,
$\gamma''\in\Gamma_{\Lambda''}$:
\begin{eqnarray*}
W_\phi(\gamma', \gamma'') \ge-\sum_{r',r''\in{\mathbb
Z}^d}\Psi(|r'-r''|_{\max}) \#(\gamma_{r'}')\, \#(\gamma_{r''}'').
\end{eqnarray*}
Here $|\cdot|_{\max}$ denotes the maximum norm on ${\mathbb R}^d$.
\end{description}
Moreover, (T) and (BB) imply
\begin{eqnarray}\label{int}
J(\beta) := \int_{{\mathbb R}^{d}} | \exp(-\beta\phi(x)) - 1 |
\,dx < \infty
\end{eqnarray}
for all $\beta \ge 0$ ($dx$ denotes the Lebesgue measure on ${\mathbb R}^d$). 
The property \eqref{int} is also called
{\it integrability} (I) or {\it regularity}.

On $(\Gamma^{\scriptscriptstyle{(N)}}_{\scriptscriptstyle{\Lambda}}, {\mathcal
B}(\Gamma^{\scriptscriptstyle{(N)}}_{\scriptscriptstyle{\Lambda}}))$ we consider the canonical
$N$-particle Gibbs measures $\mu^{\scriptscriptstyle{(N)}}_{\scriptscriptstyle{\Lambda}}$ in $\Lambda$
with empty boundary condition:
\begin{eqnarray*}
\mu^{\scriptscriptstyle{(N)}}_{\scriptscriptstyle{\Lambda}} :=
\frac{1}{Z^{\scriptscriptstyle{(N)}}_{\scriptscriptstyle{\Lambda}}}
\exp \left(-\beta  E_{\phi}\right)
dx^{\scriptscriptstyle{(N)}}_{\scriptscriptstyle{\Lambda}},
\end{eqnarray*}
where
\begin{eqnarray*}
Z^{\scriptscriptstyle{(N)}}_{\scriptscriptstyle{\Lambda}} :=
\int_{\Gamma^{\scriptscriptstyle{(N)}}_{\scriptscriptstyle{\Lambda}}}
\exp \left(-\beta E_{\phi}\right)
dx^{\scriptscriptstyle{(N)}}_{\scriptscriptstyle{\Lambda}}
\end{eqnarray*}
is the canonical partition function of $N$ particles in $\Lambda$.
The constant $\beta \ge 0$ is the inverse temperature.

For $1 \le n \le N$ the $n$-order correlation function
corresponding to $\mu^{\scriptscriptstyle{(N)}}_{\scriptscriptstyle{\Lambda}}$ is defined by
\begin{multline*}
k_{\scriptscriptstyle{\Lambda}}^{\scriptscriptstyle{(n,N)}}(x_1,
\ldots, x_n) := \frac{N \cdot \ldots \cdot
(N-n+1)}{Z^{\scriptscriptstyle{(N)}}_{\scriptscriptstyle{\Lambda}}}
\int_{\Lambda^{N-n}} \exp \Big(- \beta (E_\phi(X \cup Y) \Big)
dy_{\scriptscriptstyle{\Lambda}}^{\otimes(N-n)},
\end{multline*}
where $X = \{x_1, \ldots, x_n\}$ and $Y = \{y_1, \ldots,
y_{N-n}\}$. Furthermore, we define
\begin{eqnarray*}
k_{\scriptscriptstyle{\Lambda}}^{\scriptscriptstyle{(0,N)}} := 1
\quad \mbox{and} \quad
k_{\scriptscriptstyle{\Lambda}}^{\scriptscriptstyle{(n,N)}} := 0
\quad \mbox{for} \quad n
> N.
\end{eqnarray*}

Let $f^{(n)}:\Lambda^n \to [0,\infty]$ be a symmetric measurable
function, $1 \le n \le N$, then
\begin{align}\label{corre}
&\int_{\Gamma^{\scriptscriptstyle{(N)}}_{\scriptscriptstyle{\Lambda}}}
\sum_{\{x_1,\dots,x_n\}\subset\gamma}
f^{(n)}(x_1,\dots,x_n)\,d\mu^{\scriptscriptstyle{(N)}}_{\scriptscriptstyle{\Lambda}}(\gamma)\notag\\&\qquad
= \frac{1}{n!}\, \int_{\Lambda^n} f^{(n)}(x_1,\dots,x_n)
k_{\scriptscriptstyle{\Lambda}}^{\scriptscriptstyle{(n,N)}}(x_1,\dots,
x_n) \,dx_{\scriptscriptstyle{\Lambda}}^{\otimes n}.
\end{align}

\begin{definition}
{\rm Let $(\Lambda_N)_{N \in {\mathbb N}}$ be a sequence of
bounded Borel measurable subsets of ${\mathbb R}^d$ with
$|\Lambda_N| >0$ which exhausts ${\mathbb R}^d$, i.e., for each
bounded $\Lambda \subset {\mathbb R}^d$ there exists
$N_{\scriptscriptstyle{\Lambda}} \in {\mathbb N}$ such that
$\Lambda \subset \Lambda_N$ for all $N \ge
N_{\scriptscriptstyle{\Lambda}}$.
We denote by $|A|$ the Lebesgue measure of a Borel measurable set
$A \subset {\mathbb R}^d$. We say that $(\Lambda_N)_{N \in
{\mathbb N}}$ has an N/V-limit, if
\begin{eqnarray*}
\rho := \lim_{N \to \infty} \frac{N}{|\Lambda_N|}
\end{eqnarray*}
exists in $(0, \infty)$. In this case we call $(\Lambda_N)_{N \in
{\mathbb N}}$ a sequence of volumes corresponding to the density
$\rho > 0$. Sometimes we use the notation
$v_{\scriptscriptstyle{N}} := N/|\Lambda_N|$.}
\end{definition}

\begin{theorem}\label{ruelle}
Suppose that the conditions {\rm (RP), (T)}, and {\rm (BB)} are
satisfied and let $(\Lambda_N)_{N \in {\mathbb N}}$ be a sequence
of volumes corresponding to the density $\rho > 0$. Then for large
enough $N_0 \in {\mathbb N}$ (such that $|\Lambda_N|$ is larger
than a critical volume for all $N \ge N_0$) there exists $\xi <
\infty$ such that
\begin{eqnarray}\label{Ruellebo}
k_{\scriptscriptstyle{\Lambda_N}}^{\scriptscriptstyle{(n,N)}}(x_1,\dots,x_n)
\le \xi^n \quad \mbox{for all} \quad N \ge N_0, \, n \in {\mathbb
N}, \, (x_1,\dots,x_n) \in (\Lambda_N)^n
\end{eqnarray}
(Ruelle bound). Moreover, for $n \ge 2$ there exists $\zeta <
\infty$ such that
\begin{multline}\label{expRuellebo}
k_{\scriptscriptstyle{\Lambda_N}}^{\scriptscriptstyle{(n,N)}}(x_1,\dots,x_n)
\le \exp\Big( - \frac{2}{n}\sum_{1 \le i < j \le n} \beta
\phi(x_i-x_j) \Big) \zeta^n \\ \mbox{for all} \quad N \ge N_0, \,
n \ge 2, \, (x_1,\dots,x_n) \in (\Lambda_N)^n
\end{multline}
(improved Ruelle bound).
\end{theorem}
{\bf Proof:} The Ruelle bound for grand canonical correlation
functions is derived in \cite[Prop.~2.6]{Rue70}. Here we adapt
that proof to canonical correlations functions. For this,
additionally, we need the following estimates for canonical
partition functions provided in \cite[Lem.~3']{DoMi67}: For
$|\Lambda_N|$ large enough there exists a constant $C_1 < \infty$
such that
\begin{eqnarray}\label{domibound}
\frac{Z^{\scriptscriptstyle{(n-1)}}_{\scriptscriptstyle{\Lambda_N}}}
{Z^{\scriptscriptstyle{(n)}}_{\scriptscriptstyle{\Lambda_N}}} \le
C_1 \frac{1}{|\Lambda_N|} \quad \mbox{for all} \quad 1 \le n \le N.
\end{eqnarray}
Note that because $(\Lambda_N)_{N \in {\mathbb N}}$ has an
N/V-limit, there exists $C_2 < \infty$ such that
\begin{eqnarray}\label{seqbound}
v_{\scriptscriptstyle{N}} \le C_2 \quad \mbox{for all} \quad N \in
{\mathbb N}.
\end{eqnarray}

Now we need to introduce some notation from \cite{Rue70}. Let
$(l_j)_{j \in \mathbb N}$ be an increasing sequence in ${\mathbb
N}$. We define
\begin{eqnarray*}
[j] := \{r \in {\mathbb Z}^d \,|\, |r|_{\max} \le l_j \}, \quad
V_j := \sum_{r \in [j]} Q(r).
\end{eqnarray*}
Furthermore, let $\psi$ be an increasing function on ${\mathbb N}$
such that
\begin{eqnarray*}
\psi \ge 1, \quad \lim_{j \to \infty} \psi(j) = \infty, \quad
\mbox{and} \quad \sum_{r \in {\mathbb Z}^d} \psi(|r|_{\max})
\Psi(|r|_{\max}) < \infty.
\end{eqnarray*}
Define $\psi_j := \psi(l_j)$ and let $P > 0$. Then for each $X
\cup Y$, $X = \{x_1, \ldots, x_n\}, Y = \{y_1, \ldots, y_{N-n}\}$
either
\begin{eqnarray}\label{eqall}
\sum_{r \in [j]} \#(\{X \cup Y\}_r)^2 \le \psi_j |V_j|
\end{eqnarray}
for all $j \ge P$ or there exists a largest $q \ge P$ such that
\begin{eqnarray}\label{eqsome}
\sum_{r \in [q]} \#(\{X \cup Y\}_r)^2 \ge \psi_q |V_q|.
\end{eqnarray}
Now let $P$, the sequence $(l_j)_{j \in \mathbb N}$ and function
$\psi$ be chosen as in \cite[Sect.~2]{Rue70}. Then there exists
$C_3 < \infty$ such
\begin{eqnarray}\label{eqallest}
-W_\phi(\{x_1\}, X \setminus \{x_1\} \cup Y) \le C_3
\end{eqnarray}
for all $X \cup Y$ fulfilling \eqref{eqall}, see
\cite[Eq.~(2.29)]{Rue70}. On the other hand, for all $X \cup Y$
fulfilling \eqref{eqsome} there exists $C_4 > 0$ such that
\begin{eqnarray}\label{eqsomeest}
- E_\phi(\{X \cup Y\}_{V_{q+1}}) -W_\phi (\{X \cup Y\}_{V_{q+1}},
\{X \cup Y\}_{V^c_{q+1}}) \nonumber \\
\le - \frac{D}{4} \sum_{r \in [q+1]} \#(\{X \cup Y\}_r)^2 -
C_4 \psi_{q+1} |V_{q+1}|,
\end{eqnarray}
where the constant $D$ is as in (SS), see \cite[Prop.~2.5]{Rue70}.

We prove the assertion by induction. Let us fix $X = \{x_1,
\ldots, x_n \}, \, n \ge 1,$ and choose the coordinates of
${\mathbb Z}^d$ such that $x_1 \in Q(0)$. Let $S_0 \subset
\Lambda^{N-n}_N$ such that $X \cup Y$ fulfills \eqref{eqall} for
all $(y_1, \ldots, y_{N-n}) \in S_0$ and $S_q \subset
\Lambda^{N-n}_N$ such that $X \cup Y$ fulfills \eqref{eqsome} for
all $(y_1, \ldots, y_{N-n}) \in S_p$. Furthermore, we define
\begin{eqnarray*}
C_5 := \max \bigg\{ \bigg(\exp(\beta C_3)C_1 + \sum_{q
\ge P} \exp \Big(-(\beta C_4 \psi_{q+1} -1) |V_{q+1}|
\Big)\bigg)C_2, C_2, 1 \bigg\}
\end{eqnarray*}
and assume $\xi \ge C_5$.

Now \eqref{eqallest} together with \eqref{domibound} and
\eqref{seqbound} implies:
\begin{multline}\label{s0}
\frac{N \cdot \ldots \cdot
(N-n+1)}{Z^{\scriptscriptstyle{(N)}}_{\scriptscriptstyle{\Lambda_N}}}
\int_{S_0} \exp\big(- \beta E_\phi(X \cup Y) \big)
dy_{\Lambda_N}^{\otimes(N-n)}
\\
\le \frac{N \cdot \ldots \cdot
(N-n+1)}{Z^{\scriptscriptstyle{(N)}}_{\scriptscriptstyle{\Lambda_N}}}
\exp(\beta C_3) \int_{S_0} \exp\big(- \beta E_\phi(X \setminus
\{x_1\} \cup Y) \big) dy_{\Lambda_N}^{\otimes(N-n)}
\\
\le \exp(\beta C_3) C_1C_2
k_{\scriptscriptstyle{\Lambda_N}}^{\scriptscriptstyle{(n-1,N-1)}}
(x_2, \ldots, x_n) \le \exp(\beta C_3) C_1C_2 \xi^{n-1}.
\end{multline}

In turn, \eqref{eqsomeest} together with \eqref{domibound} and
\eqref{seqbound} yields:
\begin{multline}\label{sq}
\frac{N \cdot \ldots \cdot
(N-n+1)}{Z^{\scriptscriptstyle{(N)}}_{\scriptscriptstyle{\Lambda_N}}}\int_{S_q}
\exp\big(- \beta E_\phi(X \cup Y) \big)
dy_{\Lambda_N}^{\otimes(N-n)}
\\ \le \frac{N \cdot \ldots \cdot (N-n+1)}{Z^{\scriptscriptstyle{(N)}}_{\scriptscriptstyle{\Lambda_N}}}
\exp\left(- \beta \frac{D}{4} \sum_{r \in
V_{q+1}} \#(\{X \cup Y\}_r)^2 - \beta C_4 \psi_{q+1} |V_{q+1}| \right) \\
\times \int_{\Lambda_N^{N-n}} \exp\big( - \beta E_\phi(\{X \cup
Y\}_{V^c_{q+1}}))\big) dy_{\Lambda_N}^{\otimes(N-n)}
\\
\le \frac{N \cdot \ldots \cdot
(N-n+1)}{Z^{\scriptscriptstyle{(N)}}_{\scriptscriptstyle{\Lambda_N}}}
\exp\left(- \beta \frac{D}{4} \sum_{r \in V_{q+1}} \#(\{X \cup
Y\}_r)^2 - \beta C_4
\psi_{q+1} |V_{q+1}| \right) \\
|\Lambda_N|^{\#(Y_{V_{q+1}})}
\int_{\Lambda_N^{N-n-\#(Y_{V_{q+1}})}} \exp\big( - \beta
E_\phi(\{X \cup Y\}_{V^c_{q+1}}))\big)
dy_{\Lambda_N}^{\otimes(N-n-\#(Y_{V_{q+1}}))} \\
\le \exp\left(- \beta \frac{D}{4} \sum_{r \in V_{q+1}} \#(\{X \cup
Y\}_r)^2
- \beta C_4 \psi_{q+1} |V_{q+1}| \right) \\
\times (C_1)^{\#(Y_{V_{q+1}})}
(C_1|v_{\scriptscriptstyle{N}}|)^{\#(X_{V_{q+1}})}
k_{\scriptscriptstyle{\Lambda_N}}^{{(n -
\#(X_{V_{q+1}}),\#(Y_{V^c_{q+1}}))}}(x_{i_1}, \ldots,
x_{i_{n-\#(X_{V_{q+1}})}})
\\
\le \exp\bigg(- \beta \frac{D}{4} \sum_{r \in V_{q+1}} \#(\{X \cup
Y\}_r)^2  + \ln(C_1) \sum_{r \in V_{q+1}} \#(\{X \cup Y\}_r)\\
- \beta C_4 \psi_{q+1} |V_{q+1}| \bigg) C_2 \xi^{n-1} \le
\exp\left(-(\beta C_4\psi_{q+1} -1) |V_{q+1}| \right) C_2
\xi^{n-1},
\end{multline}
where we used that
\begin{eqnarray*}
- \beta \frac{D}{4} \#(\{X \cup Y\}_r)^2 + \ln(C_1) \#(\{X
\cup Y\}_r) \le 1.
\end{eqnarray*}
Finally, summing up \eqref{s0} and \eqref{sq} we get
\begin{multline*}
k_{\scriptscriptstyle{\Lambda_N}}^{\scriptscriptstyle{(n,N)}}
(x_1, \ldots, x_n) \le \bigg(\exp(\beta C_3)C_1 + \sum_{q \ge P}
\exp\left(-(\beta C_4\psi_{q+1} -1) |V_{q+1}| \right) \bigg) C_2
\xi^{n-1} \le \xi^{n}.
\end{multline*}

The canonical correlation functions fulfill the following
Kirkwood--Salsburg type equations:
\begin{multline*}
k_{\scriptscriptstyle{\Lambda_N}}^{\scriptscriptstyle{(n,N)}}
(x_1, \ldots, x_n) = N
\frac{Z^{\scriptscriptstyle{(N-1)}}_{\scriptscriptstyle{\Lambda_N}}}
{Z^{\scriptscriptstyle{(N)}}_{\scriptscriptstyle{\Lambda_N}}}
\exp\Big(-\sum_{2 \le i \le n} \beta \phi(x_1-x_j) \Big)
\Bigg(k_{\scriptscriptstyle{\Lambda_N}}^{\scriptscriptstyle{(n-1,N-1)}}(x_2,
\ldots, x_n) \\ + \sum_{k=1}^{N-n} \frac{1}{k!} \int_{\Lambda_N^k}
k_{\scriptscriptstyle{\Lambda_N}}^{\scriptscriptstyle{(n+k-1,N-1)}}
(x_2, \ldots, x_n, y_1, \ldots, y_k) \prod_{i=1}^k (\exp(-\beta
\phi(x_1-y_i))-1)\,dy_{\scriptscriptstyle{\Lambda}}^{\otimes k}
\Bigg),
\end{multline*}
see e.g.~\cite[Eq.~(38.16)]{Hill56}. We set
\begin{eqnarray*}
\zeta := \max \Big\{ C_1C_2 \exp(\xi I), \xi \Big\}.
\end{eqnarray*}
Then \eqref{domibound} and \eqref{seqbound} together with
the Ruelle bound \eqref{Ruellebo} yield
\begin{multline*}
k_{\scriptscriptstyle{\Lambda_N}}^{\scriptscriptstyle{(n,N)}}(x_1,
\ldots, x_n) \le \exp\Big(-\sum_{2 \le i \le n} \beta
\phi(x_1-x_j) \Big) C_1 C_2
\Big( \xi^{n-1} + \sum_{k=1}^{N-n} \frac{1}{k!} \xi^{n+k-1} I^k\Big) \\
\le \exp\Big(-\sum_{2 \le i \le n} \beta \phi(x_1-x_j) \Big) \xi^{n-1}
C_1 C_2 \exp(\xi I)
\le \exp\Big(-\sum_{2 \le i \le n} \beta \phi(x_1-x_j) \Big) \zeta^{n}.
\end{multline*}
Finally, symmetry of the correlation functions gives \eqref{expRuellebo}.
\hfill$\blacksquare$

\section{n-particle stochastic dynamics in finite volume}\label{s3}

Let $\Lambda \subset {\mathbb R}^d$ such that $\Lambda^N \subset
{\mathbb R}^{N \cdot d}$ is the closure of an open, relatively
compact set, having boundary $\partial(\Lambda^N)$ of Lebesgue
measure zero. Our aim is to construct an $N$-particle diffusion
process $({\bf X}(t))_{t \ge 0}$ in
$\Gamma_{\scriptscriptstyle{\Lambda}}^{\scriptscriptstyle{(N)}}$
solving weakly the following $N$-system of stochastic differential
equations before hitting $\partial
(\Gamma_{\scriptscriptstyle{\Lambda}}^{\scriptscriptstyle{(N)}})$:
\begin{align}
dx(t) & = - \beta \!\! \sum_{\stackunder{y(t) \, \neq \, x(t)}
{y(t) \in \,\, {\bf X}(t)}}
\nabla \phi(x(t) - y(t)) \,dt + \sqrt{2} \,dB^{x_0}(t), \nonumber \\
& \mbox{with reflecting boundary condition}, \label{npsd}
\end{align}
for sufficiently many initial conditions $\gamma_0 \in
\Gamma_{\scriptscriptstyle{\Lambda}}^{\scriptscriptstyle{(N)}}$.
Here $x(t) \in {\bf X}(t) \in
\Gamma_{\scriptscriptstyle{\Lambda}}^{\scriptscriptstyle{(N)}}$
and $(B^{x_0})_{x_0 \in \gamma_0}$ are $N$ independent Brownian
motions starting in $x_0$. Existence of a solution to \eqref{npsd}
was shown in \cite{FaGr04a} by using Dirichlet form techniques.
Here we briefly summarize their construction.

First we have to introduce an additional condition:
\begin{description}
\item[(D)] ({\it Differentiability}) The function $\exp(-\beta\phi)$ is
weakly differentiable on ${\mathbb R}^d$, $\phi$ is continuously
differentiable on ${\mathbb R}^d \backslash {\{} 0 {\}}$ and the
gradient $\nabla \phi$, considered as a $dx$-a.e.~defined function
on ${\mathbb R}^d$, satisfies
\begin{eqnarray*}
\nabla \phi \in L^1({\mathbb R}^d, \exp(-\beta\phi)dx) \cap
L^2({\mathbb R}^d, \exp(-\beta\phi)dx) \cap
L^3({\mathbb R}^d, \exp(-\beta\phi)dx),
\end{eqnarray*}
where $\beta > 0$ is the inverse temperature. 
Furthermore, we assume $\phi$ to be such that the function $\Phi$
in (RP) can be chosen differentiable and $\Phi'\exp(-a\Phi)$ a
bounded function on $(0, \infty)$ for all $a > 0$.
\end{description}
Note that, for many typical potentials in Statistical Physics, we
have $\phi \in C^{\infty}({\mathbb R}^d\backslash {\{} 0 {\}})$.
For such ``outside the origin regular'' potentials, condition (D)
nevertheless does not exclude a singularity at the point $0 \in
{\mathbb R}^d$. The last assumption on $\phi$, ensuring a suitable
choice of $\Phi$, is no restriction from the physical point of
view. E.g., potentials, diverging faster than $\Phi(t) =
t^{-d-\epsilon}$, $\epsilon > 0$, at the origin, are admissible.

On $\Lambda^N$ consider the measure
\begin{eqnarray*}
\mu_{\scriptscriptstyle{\Lambda, N}} :=
\frac{1}{Z^{\scriptscriptstyle{(N)}}_{\scriptscriptstyle{\Lambda}}}
\exp \Big(-\beta  \sum_{1 \le i<j \le N} \phi(x_i-x_j)\Big)
dx^{\otimes N}_{\scriptscriptstyle{\Lambda}}.
\end{eqnarray*}
Note that then $\mbox{sym}^{\scriptscriptstyle{(N)}}: \Lambda^N
\to
\Gamma^{\scriptscriptstyle{(N)}}_{\scriptscriptstyle{\Lambda}}$ is
$\mu_{\scriptscriptstyle{\Lambda, N}}$-a.e.~defined (since the
diagonals have $\mu_{\scriptscriptstyle{\Lambda, N}}$-measure
zero) and that
$\mu_{\scriptscriptstyle{\Lambda}}^{\scriptscriptstyle{(N)}} =
\mu_{\scriptscriptstyle{\Lambda, N}} \circ
(\mbox{sym}^{\scriptscriptstyle{(N)}})^{-1}$. Denote by $\nabla_i$
the gradient on ${\mathbb R}^d$ w.r.t. the variable $x_i$. Then
\begin{eqnarray}\label{diri}
{\mathcal E}_{\scriptscriptstyle{\Lambda, N}}(F,G) := \sum_{1 \le
i \le N} \int_{\Lambda^N} (\nabla_iF, \nabla_iG)_{{\mathbb R}^d}
\,d\mu_{\scriptscriptstyle{\Lambda, N}}, \quad F,G \in D,
\end{eqnarray}
defines a bilinear from on
\begin{eqnarray}\label{diridomain}
D = \Big\{ F \in  C(\Lambda^N) \, \Big| \, \nabla_i F
\,\,\mbox{locally} \,\, dx\mbox{-integrable on} \,
\mathring{\Lambda^N}, \nabla_i F \in
L^2(\mu_{\scriptscriptstyle{\Lambda, N}}) \Big\}.
\end{eqnarray}
Here $(\cdot,\cdot)_{{\mathbb R}^d}$ denotes the scalar product in
${\mathbb R}^d$ inducing the Euclidean norm, $\mathring{A}$ the
open kernel of the set $A$, and the gradient $\nabla_iF$ is meant
in the distributional sense on $\mathring{\Lambda^N}$. $({\mathcal
E}_{\scriptscriptstyle{\Lambda, N}}, D)$ is a densely defined,
positive definite, symmetric bilinear form on
$L^2(\mu_{\scriptscriptstyle{\Lambda, N}})$.

In \cite[Prop.~5.3]{FaGr04a} it was shown that $({\mathcal
E}_{\scriptscriptstyle{\Lambda, N}}, D)$ is closable and its
closure $({\mathcal E}_{\scriptscriptstyle{\Lambda, N}},$
$D({\mathcal E}_{\scriptscriptstyle{\Lambda, N}}))$ is a
conservative, local, quasi-regular Dirichlet form. Thus, there exists a
corresponding self-adjoint generator
$(H_{\scriptscriptstyle{\Lambda, N}},
D(H_{\scriptscriptstyle{\Lambda, N}}))$ (Friedrichs extension),
i.e., $D(H_{\scriptscriptstyle{\Lambda, N}}) \subset D({\mathcal
E}_{\scriptscriptstyle{\Lambda, N}})$ and
\begin{eqnarray*}
{\mathcal E}_{\scriptscriptstyle{\Lambda, N}}(F,G) =
\int_{\Lambda^{N}} H_{\scriptscriptstyle{\Lambda, N}} F G \,d \mu_{\scriptscriptstyle{\Lambda, N}},
\qquad F \in D(H_{\scriptscriptstyle{\Lambda, N}}), \,
G \in D({\mathcal E}_{\scriptscriptstyle{\Lambda, N}}).
\end{eqnarray*}

In order to solve \eqref{npsd}, however, we are rather interested
in the image Dirichlet form under
$\mbox{sym}^{\scriptscriptstyle{(N)}}$. Define an isometry
$(\mbox{sym}^{\scriptscriptstyle{(N)}})^*:
L^2(\Gamma_{\scriptscriptstyle{\Lambda}}^{\scriptscriptstyle{(N)}},
\mu^{\scriptscriptstyle{(N)}}_{\scriptscriptstyle{\Lambda}}) \to
L^2(\Lambda^N, \mu_{\scriptscriptstyle{\Lambda, N}})$ by setting
$(\mbox{sym}^{\scriptscriptstyle{(N)}})^*F$ to be the
$\mu_{\scriptscriptstyle{\Lambda, N}}$-class represented by
$\tilde{F} \circ \mbox{sym}^{\scriptscriptstyle{(N)}}$ on
$\widetilde{\Lambda^N}$ for any
$\mu^{\scriptscriptstyle{(N)}}_{\scriptscriptstyle{\Lambda}}$-version
$\tilde{F}$ of $F \in
L^2(\mu^{\scriptscriptstyle{(N)}}_{\scriptscriptstyle{\Lambda}})$.
Note that the subspace
\begin{eqnarray*}
L_{\rm sym}^2(\mu_{\scriptscriptstyle{\Lambda, N}})
:= (\mbox{sym}^{\scriptscriptstyle{(N)}})^*
(L^2(\mu^{\scriptscriptstyle{(N)}}_{\scriptscriptstyle{\Lambda}}))
\subset L^2(\mu_{\scriptscriptstyle{\Lambda, N}})
\end{eqnarray*}
is the closed subspace of symmetric functions from
$L^2(\mu_{\scriptscriptstyle{\Lambda, N}})$. Using this mapping
one can define a bilinear form $({\mathcal
E}^{\scriptscriptstyle{(N)}}_{\scriptscriptstyle{\Lambda}},
D({\mathcal
E}^{\scriptscriptstyle{(N)}}_{\scriptscriptstyle{\Lambda}}))$ as
the image bilinear form of $({\mathcal
E}_{\scriptscriptstyle{\Lambda, N}}, D({\mathcal
E}_{\scriptscriptstyle{\Lambda, N}}))$ under
$\mbox{sym}^{\scriptscriptstyle{(N)}}$:
\begin{align}\label{dfconfi}
D({\mathcal
E}^{\scriptscriptstyle{(N)}}_{\scriptscriptstyle{\Lambda}}) & := 
\{ F \in L^2(\mu^{\scriptscriptstyle{(N)}}_{\scriptscriptstyle{\Lambda}})   
\,|\, ({\rm sym}^{\scriptscriptstyle{(N)}})^*F \in 
D({\mathcal E}_{\scriptscriptstyle{\Lambda, N}}) \},
\nonumber \\ {\mathcal
E}^{\scriptscriptstyle{(N)}}_{\scriptscriptstyle{\Lambda}}(F,G) & :=
{\mathcal E}_{\scriptscriptstyle{\Lambda, N}}(({\rm
sym}^{\scriptscriptstyle{(N)}})^*F, ({\rm
sym}^{\scriptscriptstyle{(N)}})^*G), \qquad F, G \in D({\mathcal
E}^{\scriptscriptstyle{(N)}}_{\scriptscriptstyle{\Lambda}}).
\end{align}
Also $({\mathcal
E}^{\scriptscriptstyle{(N)}}_{\scriptscriptstyle{\Lambda}},
D({\mathcal
E}^{\scriptscriptstyle{(N)}}_{\scriptscriptstyle{\Lambda}}))$ is a
conservative, local, symmetric Dirichlet form. Its generator is
given by
\begin{align}\label{genconfi}
H^{\scriptscriptstyle{(N)}}_{\scriptscriptstyle{\Lambda}} &=
(({\rm sym}^{\scriptscriptstyle{(N)}})^*)^{-1}
 \circ H_{\scriptscriptstyle{\Lambda, N}} \circ ({\rm
sym}^{\scriptscriptstyle{(N)}})^*, \nonumber \\
D(H^{\scriptscriptstyle{(N)}}_{\scriptscriptstyle{\Lambda}}) &= 
\{ F \in L^2(\mu^{\scriptscriptstyle{(N)}}_{\scriptscriptstyle{\Lambda}})   
\,|\, ({\rm sym}^{\scriptscriptstyle{(N)}})^*F \in 
D(H_{\scriptscriptstyle{\Lambda, N}}) \}.
\end{align}
Of course, $(H_{\scriptscriptstyle{\Lambda, N}},
D(H_{\scriptscriptstyle{\Lambda, N}}))$ generates a strongly
continuous contraction semi-group
\begin{eqnarray*}
T^{\scriptscriptstyle{(N)}}_{\scriptscriptstyle{\Lambda}}(t) :=
\exp(-tH^{\scriptscriptstyle{(N)}}_{\scriptscriptstyle{\Lambda}}),
\quad t \ge 0.
\end{eqnarray*}

For repulsive potentials satisfying (D) in \cite{FaGr04a} it was
shown that $Dg := \Lambda^N \setminus \widetilde{\Lambda^N}$
has ${\mathcal
E}_{\scriptscriptstyle{\Lambda, N}}$-capacity
zero. Thus, $({\mathcal
E}^{\scriptscriptstyle{(N)}}_{\scriptscriptstyle{\Lambda}},
D({\mathcal
E}^{\scriptscriptstyle{(N)}}_{\scriptscriptstyle{\Lambda}}))$ is
obviously also quasi-regular and by \cite[Chap.~IV, Sect.~3]{MaRo92}
we have the following theorem:

\begin{theorem}\label{exnpsd}
Suppose that conditions {\rm (RP)}, {\rm (D)} are satisfied, $N
\in {\mathbb N}$ and $\Lambda \subset {\mathbb R}^d$ such that $\Lambda^N
\subset {\mathbb R}^{N \cdot d}$ is the closure of an open, relatively
compact set with boundary
$\partial(\Lambda^N)$ of Lebesgue measure zero. Then:

\noindent
{\rm (i)} There exists a conservative diffusion process
(i.e., a conservative strong Markov process with continuous sample paths)
\begin{eqnarray*}
{\bf M}_{\scriptscriptstyle{\Lambda}}^{\scriptscriptstyle{(N)}} =
\big({\bf
\Omega}_{\scriptscriptstyle{\Lambda}}^{\scriptscriptstyle{(N)}},
{\bf F}_{\scriptscriptstyle{\Lambda}}^{\scriptscriptstyle{(N)}},
({\bf
F}_{\scriptscriptstyle{\Lambda}}^{\scriptscriptstyle{(N)}}(t))_{t
\ge 0}, ({\bf
\Theta}_{\scriptscriptstyle{\Lambda}}^{\scriptscriptstyle{(N)}}(t))_{t
\ge 0}, ({\bf X}(t))_{t \ge 0}, ({\bf
P}_{\scriptscriptstyle{\Lambda}}^{\scriptscriptstyle{(N)}}(x))_{x
\in
\Gamma^{\scriptscriptstyle{(N)}}_{\scriptscriptstyle{\Lambda}}}\big)
\end{eqnarray*}
on
$\Gamma^{\scriptscriptstyle{(N)}}_{\scriptscriptstyle{\Lambda}}$
which is properly associated with $({\mathcal
E}_{\scriptscriptstyle{\Lambda}}^{\scriptscriptstyle{(N)}},
D({\mathcal
E}_{\scriptscriptstyle{\Lambda}}^{\scriptscriptstyle{(N)}}))$,
i.e., for all
($\mu_{\scriptscriptstyle{\Lambda}}^{\scriptscriptstyle{(N)}}$-versions
of) $F \in
L^2(\Gamma^{\scriptscriptstyle{(N)}}_{\scriptscriptstyle{\Lambda}},
\mu_{\scriptscriptstyle{\Lambda}}^{\scriptscriptstyle{(N)}})$ and
all $t > 0$ the function
\begin{eqnarray*}
x \mapsto \int_{{\bf
\Omega}_{\scriptscriptstyle{\Lambda}}^{\scriptscriptstyle{(N)}}}
F({\bf X}(t)) \,d{\bf
P}_{\scriptscriptstyle{\Lambda}}^{\scriptscriptstyle{(N)}}(x),
\qquad x \in
\Gamma^{\scriptscriptstyle{(N)}}_{\scriptscriptstyle{\Lambda}},
\end{eqnarray*}
is a ${\mathcal
E}^{\scriptscriptstyle{(N)}}_{\scriptscriptstyle{\Lambda}}$-quasi-continuous
version of
$T_{\scriptscriptstyle{\Lambda}}^{\scriptscriptstyle{(N)}}(t)F$.
${\bf M}_{\scriptscriptstyle{\Lambda}}^{\scriptscriptstyle{(N)}}$
is up to
$\mu_{\scriptscriptstyle{\Lambda}}^{\scriptscriptstyle{(N)}}$-equivalence
unique. In particular, ${\bf
M}_{\scriptscriptstyle{\Lambda}}^{\scriptscriptstyle{(N)}}$ is
$\mu_{\scriptscriptstyle{\Lambda}}^{\scriptscriptstyle{(N)}}$-symmetric
(i.e., $\int G \,
T_{\scriptscriptstyle{\Lambda}}^{\scriptscriptstyle{(N)}}(t) F \,
d\mu_{\scriptscriptstyle{\Lambda}}^{\scriptscriptstyle{(N)}} =
\int F \,
T_{\scriptscriptstyle{\Lambda}}^{\scriptscriptstyle{(N)}}(t) G
\,d\mu_{\scriptscriptstyle{\Lambda}}^{\scriptscriptstyle{(N)}}$
for all $F,G:
\Gamma^{\scriptscriptstyle{(N)}}_{\scriptscriptstyle{\Lambda}}
 \to [0, \infty)$ measurable) and  has
$\mu_{\scriptscriptstyle{\Lambda}}^{\scriptscriptstyle{(N)}}$ as
an invariant measure.

\noindent {\rm (ii)} The diffusion process ${\bf
M}_{\scriptscriptstyle{\Lambda}}^{\scriptscriptstyle{(N)}}$ is up
to
$\mu_{\scriptscriptstyle{\Lambda}}^{\scriptscriptstyle{(N)}}$-equivalence
the unique diffusion process having
$\mu_{\scriptscriptstyle{\Lambda}}^{\scriptscriptstyle{(N)}}$ as
symmetrizing measure and solving the martingale problem for
$(-H_{\scriptscriptstyle{\Lambda}}^{\scriptscriptstyle{(N)}},$
$D(H_{\scriptscriptstyle{\Lambda}}^{\scriptscriptstyle{(N)}}))$,
in the sense that for all $G \in
D(H_{\scriptscriptstyle{\Lambda}}^{\scriptscriptstyle{(N)}})$
\begin{eqnarray*}
G({\bf X}(t)) - G({\bf X}(0)) + \int_0^t
H_{\scriptscriptstyle{\Lambda}}^{\scriptscriptstyle{(N)}}G({\bf
X}(s))\,ds, \qquad t \ge 0,
\end{eqnarray*}
is an ${\bf
F}_{\scriptscriptstyle{\Lambda}}^{\scriptscriptstyle{(N)}}(t)$-martingale
under ${\bf
P}_{\scriptscriptstyle{\Lambda}}^{\scriptscriptstyle{(N)}}(x)$
(hence starting in $x$) for ${\mathcal
E}^{\scriptscriptstyle{(N)}}_{\scriptscriptstyle{\Lambda}}$-quasi
all $x \in
\Gamma^{\scriptscriptstyle{(N)}}_{\scriptscriptstyle{\Lambda}}$.
\end{theorem}

In the above theorem ${\bf
M}_{\scriptscriptstyle{\Lambda}}^{\scriptscriptstyle{(N)}}$ is
canonical, i.e., ${\bf
\Omega}_{\scriptscriptstyle{\Lambda}}^{\scriptscriptstyle{(N)}} =
C([0, \infty) \to
\Gamma^{\scriptscriptstyle{(N)}}_{\scriptscriptstyle{\Lambda}}),
\, {\bf X}(t)(\omega) = \omega(t), \omega \in {\bf
\Omega}_{\scriptscriptstyle{\Lambda}}^{\scriptscriptstyle{(N)}}$.
The filtration $({\bf
F}_{\scriptscriptstyle{\Lambda}}^{\scriptscriptstyle{(N)}}(t))_{t
\ge 0}$ is the natural ``minimum completed admissible
filtration'', cf.~\cite{FOT94}, Chap.~A.2, or \cite{MaRo92},
Chap.~IV, obtained from the $\sigma$-algebras $\sigma{\{}
\omega(s) \, | \, 0 \le s \le t, \omega \in {\bf
\Omega}_{\scriptscriptstyle{\Lambda}}^{\scriptscriptstyle{(N)}}
{\}}, \, t \ge 0$. ${\bf
F}_{\scriptscriptstyle{\Lambda}}^{\scriptscriptstyle{(N)}} := {\bf
F}_{\scriptscriptstyle{\Lambda}}^{\scriptscriptstyle{(N)}}(\infty)
:= \bigvee_{t \in [0, \infty)} {\bf
F}_{\scriptscriptstyle{\Lambda}}^{\scriptscriptstyle{(N)}}(t)$ is
the smallest $\sigma$-algebra containing all ${\bf
F}_{\scriptscriptstyle{\Lambda}}^{\scriptscriptstyle{(N)}}(t)$ and
$({\bf
\Theta}_{\scriptscriptstyle{\Lambda}}^{\scriptscriptstyle{(N)}}(t))_{t
\ge 0}$ are the corresponding natural time shifts. For a detailed
discussions of these objects we refer to \cite{MaRo92}.

To illustrate the relation of the process ${\bf
M}_{\scriptscriptstyle{\Lambda}}^{\scriptscriptstyle{(N)}}$ to the
stochastic differential equation \eqref{npsd} we need an explicit
representation of the generator $(H_{\scriptscriptstyle{\Lambda,
N}}, D(H_{\scriptscriptstyle{\Lambda, N}}))$, at least for some
subset of $D(H_{\scriptscriptstyle{\Lambda, N}})$. An integration by parts
yields the following representation for
$H_{\scriptscriptstyle{\Lambda, N}}$ restricted to $F \in
C_c^2(\mathring{\Lambda^N}) \subset
D(H_{\scriptscriptstyle{\Lambda, N}})$:
\begin{align}\label{repgen}
&H_{\scriptscriptstyle{\Lambda,N}}F(x) = - \sum_{i=1}^{N} \Delta_i
F(x) + \beta \!\! \sum_{1 \le i < j \le N} \nabla \phi(x_i-x_j)
(\nabla_iF(x) - \nabla_jF(x)),
\end{align}
$x \in \Lambda^N$. Furthermore, if we assume $\partial(\Lambda^N)$ to be
Lipschitz, then
\begin{eqnarray}\label{eqneu}
\Big\{F \in C^2(\Lambda^N) \,|\, \partial_{\nu} F = 0 \,\, \mbox{on} \,\,
\partial (\Lambda^N) \Big\}
\subset D(H_{\scriptscriptstyle{\Lambda, N}}),
\end{eqnarray}
where $ \partial_\nu$ denotes the normal derivative, and
representation \eqref{repgen} holds for such functions also, see
\cite[Theo.~3.2]{FaGr04a}. Note that the functions in \eqref{eqneu}
have Neumann boundary condition on $\Lambda^N$.

Let ${\mathcal F}C_{\mathrm b}^\infty({\mathcal D},\Gamma)$ be the
set of all functions on $\Gamma$ of the form
\begin{eqnarray}\label{fcbinfty}
F(\gamma)=g_F(\langle f_1,\gamma\rangle,\dots,\langle
f_n,\gamma\rangle),
\end{eqnarray}
where $n\in{\mathbb N}$, $f_1,\dots,f_n\in {\mathcal D} :=
C^{\infty}_c({\mathbb R}^d)$, and $g_F\in C^\infty_{\mathrm
b}({\mathbb R}^n)$. Here $C^{\infty}_c({\mathbb R}^d)$ denotes the
set of all infinitely differentiable functions on ${\mathbb R}^d$
with compact support and $C^\infty_{\mathrm b}({\mathbb R}^n)$
denotes the set of all infinitely differentiable functions on
${\mathbb R}^n$ which are bounded together with all their
derivatives. For $F$ as in \eqref{fcbinfty} such that
$({\rm sym}^{\scriptscriptstyle{(N)}})^*F \in C_c^2(\mathring{\Lambda^N})$,
\eqref{repgen} together with \eqref{genconfi} yields
\begin{multline}\label{explicitgenconfi}
H_{\scriptscriptstyle{\Lambda}}^{\scriptscriptstyle{(N)}}F(\gamma)
= - \sum_{i,j=1}^{N} \partial_i
\partial_j g_F(\langle f_1, \gamma \rangle, \ldots, \langle f_N,
\gamma \rangle) \langle (\nabla f_i, \nabla f_j)_{{\mathbb R}^d},
\gamma \rangle \\
- \sum_{j=1}^{N} \partial_j g_F(\langle f_1, \gamma \rangle,
\ldots, \langle f_N, \gamma \rangle) \Big{(} \langle \Delta f_j,
\gamma \rangle - \beta \!\! \sum_{\{x,y\} \subset \gamma} \nabla
\phi(x-y) (\nabla f_j(x) - \nabla f_j(y)) \Big{)},
\end{multline}
$\gamma \in
\Gamma_{\scriptscriptstyle{\Lambda}}^{\scriptscriptstyle{(N)}}$,
where $\partial_j$ denotes the partial derivative w.r.t.~the
$j$-th variable.

Now, using It\^o's formula, we find that the process ${\bf
P}_{\scriptscriptstyle{\Lambda}}^{\scriptscriptstyle{(N)}}$ solves
the stochastic differential equation \eqref{npsd} in the sense of
the associated martingale problem, see Theorem \ref{exnpsd}(ii).
We say that ${\bf
P}_{\scriptscriptstyle{\Lambda}}^{\scriptscriptstyle{(N)}}$
corresponds to reflecting boundary condition because of
the Neumann boundary
condition seen on the level of the domain of its generator, see
\eqref{eqneu}.

\section{n/v-limit of n-particle, finite volume stochastic dynamics}\label{s4}

As state space for the $N/V$-limit we consider $(\Gamma,
d_{\scriptscriptstyle{(\beta/3)}\Phi,h})$ with $\Phi$ as in
condition (RP), (D), and $h$ as in Proposition \ref{propmetric}.
$\beta$ is the inverse temperature. The laws of the equilibrium
processes
\begin{eqnarray*}
{\bf
P}_{\mu_{\scriptscriptstyle{\Lambda}}^{\scriptscriptstyle{(N)}}}
:=
\int_{\Gamma_{\scriptscriptstyle{\Lambda}}^{\scriptscriptstyle{(N)}}}
{\bf
P}_{\scriptscriptstyle{\Lambda}}^{\scriptscriptstyle{(N)}}(\gamma)
\,d\mu_{\scriptscriptstyle{\Lambda}}^{\scriptscriptstyle{(N)}}(\gamma)
\end{eqnarray*}
are probability measures on $C([0, \infty),
\Gamma_{\scriptscriptstyle{\Lambda}}^{\scriptscriptstyle{(N)}})$,
cf.~Theorem \ref{exnpsd}. Since $C([0, \infty),
\Gamma_{\scriptscriptstyle{\Lambda}}^{\scriptscriptstyle{(N)}})$
is a Borel subset of $C([0, \infty), \Gamma)$ (under the natural
embedding) with compatible measurable structures we can consider
${\bf
P}_{\mu_{\scriptscriptstyle{\Lambda}}^{\scriptscriptstyle{(N)}}}$
as a measure on $C([0, \infty), \Gamma)$. Below, $({\bf X}(t))_{t
\ge 0}$ always denotes the coordinate process in the corresponding
path space. We denote by $({\bf F}_t)_{t \ge 0}$ the natural
filtration on $C([0, \infty), \Gamma)$.

\subsection{Tightness}

\begin{theorem}\label{tight}
Suppose that the conditions {\rm (RP), (T), (BB)}, and {\rm (D)}
are satisfied and let $(\Lambda_N)_{N \in {\mathbb N}}$ be a
sequence of volumes corresponding to the density $\rho > 0$.
Furthermore, assume that the $\Lambda_N \subset {\mathbb R}^d$ are
such that $(\Lambda_N)^N \subset {\mathbb R}^{N \cdot d}$ is the
closure of an open, relatively compact set with boundary
$\partial((\Lambda_N)^N)$ of Lebesgue measure zero. Set ${\bf
P}^{\scriptscriptstyle{(N)}} := {\bf
P}_{\mu_{\scriptscriptstyle{\Lambda_N}}^{\scriptscriptstyle{(N)}}}$.
Then $({\bf P}^{\scriptscriptstyle{(N)}})_{N \in {\mathbb N}}$ is
tight on $C([0, \infty), \Gamma)$.
\end{theorem}

For a symmetric function $f: {\mathbb R}^d \times  {\mathbb R}^d \to {\mathbb R}$
we set functions $f^{[n,2]}: {\mathbb R}^{n \cdot d} \to {\mathbb R}$, $n = 2, 3, 4$, 
by
\begin{align*}
f^{[2,2]}(x_1, x_2) & := f(x_1, x_2)^2, \\
f^{[3,2]}(x_1, x_2, x_3) & := f(x_1, x_2)f(x_1, x_3) +
f(x_1, x_3)f(x_2, x_3) + f(x_2, x_3)f(x_1, x_2), \\
f^{[4,2]}(x_1, x_2, x_3, x_4) & := f(x_1, x_2)f(x_3, x_4) +
f(x_1, x_3)f(x_2, x_4) + f(x_1, x_4)f(x_2, x_3), 
\end{align*}
for $x_1, x_2, x_3, x_4 \in {\mathbb R}$, and functions 
$f^{[n,3]}: {\mathbb R}^{n \cdot d} \to {\mathbb R}$, $n = 2, \ldots, 6$, 
by
\begin{align*}
f^{[2,3]}(x_1, x_2) & := f(x_1, x_2)^3, \\
f^{[3,3]}(x_1, x_2, x_3) & :=  f(x_1, x_2)f(x_1, x_3)f(x_2, x_3) + 
f(x_1, x_2)^2\Big(f(x_1, x_3) + f(x_2, x_3) \Big) \\ 
& \! + f(x_1, x_3)^2\Big(f(x_1, x_2) + f(x_2, x_3)\Big)
+ f(x_2, x_3)^2\Big(f(x_1,x_2) + f(x_1,x_3)\Big), \\
f^{[4,3]}(x_1, \ldots, x_4) & := f(x_1, x_2)f(x_3, x_4) f(x_1, x_3) + \ldots 
+ f(x_1, x_2)^2 f(x_3, x_4) + \ldots,
\\
f^{[5,3]}(x_1, \ldots, x_5) & := f(x_1, x_2)f(x_3, x_4) f(x_5, x_1) 
+ \ldots, \\   
f^{[6,3]}(x_1, \ldots, x_6) & := f(x_1, x_2)f(x_3, x_4) f(x_5, x_6) 
+ \ldots,  
\end{align*}
for $x_1, \ldots, x_6 \in {\mathbb R}$.

\begin{lemma}\label{tight0}
{\rm Let the assumptions in Theorem \ref{tight} hold and $\beta,
\Phi, h$ be as in the metric $d_{\scriptscriptstyle{(\beta/3)}
\Phi,h}$. Set $\mu^{\scriptscriptstyle{(N)}} :=
{\mu}_{\scriptscriptstyle{\Lambda_N}}^{\scriptscriptstyle{(N)}}$.
Then
\begin{eqnarray*}
\sup_{N \in {\mathbb N}}{\mathbb
E}_{{\mu}^{\scriptscriptstyle{(N)}}}
[(S^{\scriptscriptstyle{(\beta/3)}\Phi,h})^2] < \infty.
\end{eqnarray*}
}
\end{lemma}
{\bf Proof:} Set $f(x,y) := \exp((\beta/3)\Phi(|x-y|))h(x)h(y)$,
$\{x, y\} \subset {\mathbb R}^d$. Then
\begin{multline*}
(S^{\scriptscriptstyle{(\beta/3)}\Phi,h}(\gamma))^2 =
\sum_{\{x,y,z,w\} \subset \gamma} f(x,y) f(z,w) + f(x,z)f(y,w) +
f(x,w)f(y,z) \\ + \sum_{\{x,y,z\} \subset \gamma} f(x,y) f(y,z) +
f(x,z)f(y,z) + f(x,y)f(x,z) + \sum_{\{x,y\} \subset \gamma}
f(x,y)^2.
\end{multline*}
Now \eqref{corre} together with \eqref{expRuellebo} yields for $N \ge N_0$
(as in Theorem \ref{ruelle})
\begin{multline}\label{crit1}
{\mathbb E}_{{\mu}^{\scriptscriptstyle{(N)}}}
[(S^{\scriptscriptstyle{(\beta/3)}\Phi,h})^2] 
= \sum_{n=2}^4 \frac{1}{n!} \int_{\Lambda_N^n} f^{[n,2]}(x_1, \ldots, x_n) 
k_{\scriptscriptstyle{\Lambda_N}}^{\scriptscriptstyle{(n,N)}}(x_1, \ldots, x_n) 
\,dx_{\Lambda_N}^{\otimes n} \\
\le \sum_{n=2}^4 \frac{\zeta^n}{n!}
\int_{({\mathbb R}^d)^n}  |f|^{[n,2]}(x_1, \ldots, x_n) \exp\Big( - \frac{2}{n}\sum_{1
\le i < j \le n} \beta \phi(x_i-x_j) \Big) \,dx^{\otimes n}.
\end{multline}
The integrals in \eqref{crit1} are finite due to the integrability
properties of $h$ and (RP) (note that $\exp(b\Phi(|\cdot|))
\exp(-c\phi)$ is a bounded function for all $c \ge b \ge 0$).
Therefore, \eqref{crit1} is a bound for ${\mathbb
E}_{{\mu}^{\scriptscriptstyle{(N)}}}
[(S^{\scriptscriptstyle{(\beta/3)}\Phi,h})^2]$ uniformly in $N \ge
N_0$. Of course, ${\mathbb E}_{{\mu}^{\scriptscriptstyle{(N)}}}
[(S^{\scriptscriptstyle{(\beta/3)}\Phi,h})^2]$ is finite for the
finite many $N < N_0$. Thus, the assertion is proven. \hfill
$\blacksquare$

\begin{lemma}\label{tight1}
{\rm Let the assumptions in Theorem \ref{tight} hold. Then there
exists $C_6 < \infty$ such that
\begin{eqnarray}\label{kolche}
\sup_{N \in {\mathbb N}}{\mathbb E}_{{\bf
P}^{\scriptscriptstyle{(N)}}}
\Big[d_{\scriptscriptstyle{(\beta/3)}\Phi,h}({\bf  X}(t), {\bf
X}(s))^4\Big]^{1/4} \le  C_6 \, (t-s)^{1/2}.
\end{eqnarray}
}
\end{lemma} {\bf Proof:} Recall the definition of the metric 
$d_{\scriptscriptstyle{(\beta/3)}\Phi,h}$, 
see \eqref{metric}. Since $|1-\exp(-r)| \le r$ for $r \ge 0$, by
the triangle inequality we obtain
\begin{multline}\label{kolching}
{\mathbb E}_{{\bf P}^{\scriptscriptstyle{(N)}}}
\Big[d_{\scriptscriptstyle{(\beta/3)}\Phi,h}({\bf  X}(t), {\bf
X}(s))^4\Big]^{1/4} \le \sum_{k = 1}^\infty 2^{-k} p_k {\mathbb
E}_{{\bf P}^{\scriptscriptstyle{(N)}}} [| \langle f_k, {\bf  X}(t)
\rangle - \langle f_k, {\bf X}(s) \rangle|^4]^{1/4} \\ + \sum_{k =
1}^\infty 2^{-k} q_k {\mathbb E}_{{\bf
P}^{\scriptscriptstyle{(N)}}} [|
S^{\scriptscriptstyle{(\beta/3)}\Phi,h_k}({\bf  X}(t)) -
S^{\scriptscriptstyle{(\beta/3)}\Phi,h_k}({\bf X}(s)) |^4]^{1/4}.
\end{multline}

Set $F(x) := \sum_{1 \le i \le N} f(x_i)$, $x \in (\Lambda_N)^N$,
$f \in C_c^1({\mathbb R}^d)$. By \eqref{diridomain} we know that
$F \in D({\mathcal E}_{\scriptscriptstyle{\Lambda_N, N}})$. Note
that $\langle f, \mbox{sym}^{\scriptscriptstyle{(N)}}(\cdot)
\rangle = F$ on $\widetilde{(\Lambda_N)^N}$. Thus, by
\eqref{dfconfi} $\langle f, \cdot \rangle \in D({\mathcal
E}^{\scriptscriptstyle{(N)}}_{\scriptscriptstyle{\Lambda_N}})$.
Fix $T > 0$. Below we canonically project the laws of the
equilibrium processes ${\bf P}^{\scriptscriptstyle{(N)}}$ onto
$\Omega^{\scriptscriptstyle{(N)}}_{\scriptscriptstyle{T}} := C([0,
T],
\Gamma^{\scriptscriptstyle{(N)}}_{\scriptscriptstyle{\Lambda_N}})$
without expressing this explicitly. We define the time reversal
$r_{\scriptscriptstyle{T}}(\omega) := \omega(T- \cdot), \, \omega
\in \Omega^{\scriptscriptstyle{(N)}}_{\scriptscriptstyle{T}}$.
Now, by the well-known Lyons--Zheng decomposition,
cf.~\cite{LZ88}, \cite{FOT94}, we have for all $0 \le t \le T$:
\begin{eqnarray*}
\langle f, {\bf  X}(t)\rangle - \langle f, {\bf  X}(0) \rangle =
\frac{1}{2} {\bf M}(N, f, t) + \frac{1}{2}\Big{(}{\bf M}(N, f,
T-t)(r_{\scriptscriptstyle{T}}) - {\bf  M}(N, f,
T)(r_{\scriptscriptstyle{T}}) \Big{)}
\end{eqnarray*}
${\bf P}^{\scriptscriptstyle{(N)}}$-a.e., where $({\bf M}(N, f,
t))_{0 \le t \le T}$ is a continuous $({\bf
P}^{\scriptscriptstyle{(N)}}, ({\bf
F}^{\scriptscriptstyle{(N)}}_{\scriptscriptstyle{\Lambda_N}}(t))_{0
\le t \le T})$-martingale and $({\bf M}(N, f,
t)(r_{\scriptscriptstyle{T}}))_{0 \le t \le T}$ is a continuous
$({\bf P}^{\scriptscriptstyle{(N)}},
(r_{\scriptscriptstyle{T}}^{-1} ({\bf
F}^{\scriptscriptstyle{(N)}}_{\scriptscriptstyle{\Lambda_N}}(t)))_{0
\le t \le T})$-martingale. (We note that ${\bf
P}^{\scriptscriptstyle{(N)}} \circ r_{\scriptscriptstyle{T}}^{-1}
= {\bf P}^{\scriptscriptstyle{(N)}}$ because
$(T^{\scriptscriptstyle{(N)}}_{\scriptscriptstyle{\Lambda_N}}(t))_{t
\ge 0}$ is symmetric on $L^2(\mu^{\scriptscriptstyle{(N)}})$.)
Moreover, by (\ref{diri}) the bracket of ${\bf M}(N, f)$ is given
by
\begin{align*}
<{\bf M}(N, f)>(t) & = \int_0^t \langle |\nabla f|_{{\mathbb
R}^d}^2, {\bf X}(u) \rangle \,du
\end{align*}
as e.g.~directly follows from \cite{FOT94}, Theorem 5.2.3 and
Theorem 5.1.3(i). Hence by the Burkholder--Davies--Gundy inequalities and since
${\bf P}_{\scriptscriptstyle{N}}
\circ r_{\scriptscriptstyle{T}}^{-1}
= {\bf P}_{\scriptscriptstyle{N}}$
we can find $C_7 \in (0, \infty)$
such that for all $f \in C^1_c({\mathbb R}^d), \, N \ge N_0$
(as in Theorem \ref{ruelle})$, \, 0 \le s \le t \le
T$,
\begin{multline}\label{eq1605}
{\mathbb E}_{{\bf P}^{\scriptscriptstyle{(N)}}} [| \langle f, {\bf
X}(t) \rangle - \langle f, {\bf X}(s) \rangle|^4]^{1/4} \le
\frac{1}{2} \Big( {\mathbb E}_{{\bf P}^{\scriptscriptstyle{(N)}}}
[|{\bf M}(N, f, t) - {\bf M}(N, f, s)|^4 ]^{1/4}
\\ + {\mathbb E}_{{\bf P}^{\scriptscriptstyle{(N)}}}[
|{\bf  M}(N, f, T-t)(r_{\scriptscriptstyle{T}}) - {\bf M}(N, f,
T-s)(r_{\scriptscriptstyle{T}})|^4 ]^{1/4} \Big)
\\  \le C_7\, \Bigg{(} {\mathbb
E}_{{\bf P}^{\scriptscriptstyle{(N)}}}\Big{[} \Big{(}\int_s^{t}
\langle |\nabla f|_{{\mathbb R}^d}^2, {\bf X}(u) \rangle
\,du \Big{)}^2 \Big{]}^{1/4}
\\
+  {\mathbb E}_{{\bf P}^{\scriptscriptstyle{(N)}}}\Big{[}
\Big{(}\int_{T-t}^{T-s} \langle |\nabla f|_{{\mathbb
R}^d}^2, {\bf X}(T-u) \rangle \,du \Big{)}^2 \Big{]}^{1/4} \Bigg{)} \\
\le  2 \, C_7 \,
(t-s)^{1/2} \Big{(} \int_{\Gamma_{\scriptscriptstyle{\Lambda_N}}^{\scriptscriptstyle{(N)}}}
\langle |\nabla f|_{{\mathbb R}^d}^2, \gamma \rangle^2
\,d\mu^{\scriptscriptstyle{(N)}}(\gamma) \Big{)}^{1/4} \\
\le  C_8 \, (t-s)^{1/2} \Big{(} \Big(\int_{\Lambda_N} |\nabla f|^2 \,dx_{\Lambda_N} \Big)^{2}
+ \Big(\int_{\Lambda_N} |\nabla f|^4 \,dx_{\Lambda_N} \Big) \Big)^{1/4}
\le  C_8 \, (t-s)^{1/2} I(f),
\end{multline}
where $C_8 := 2 \, C_7 \max\{ \xi^2/2, \xi\}^{1/4}$, see \eqref{corre} together with
Theorem \ref{ruelle}, and $I(f)$ is given by
\begin{eqnarray*}
I(f) = \Big{(} \Big(\int_{{\mathbb R}^d} |\nabla f|^2 \,dx \Big)^{2}
+ \Big(\int_{{\mathbb R}^d} |\nabla f|^4 \,dx \Big) \Big)^{1/4}.
\end{eqnarray*}
But then from the above derivation of \eqref{eq1605} it is clear there exists
$C_{9} < \infty$ such that
\begin{eqnarray}\label{eqc9}
{\mathbb E}_{{\bf P}^{\scriptscriptstyle{(N)}}}
[| \langle f, {\bf  X}(t) \rangle
- \langle f, {\bf X}(s) \rangle|^4]^{1/4} \le C_{9} I(f) \, (t-s)^{1/2}
\end{eqnarray}
for all $f \in C^1_c({\mathbb R}^d), \, N \in {\mathbb N}$, $0 \le s \le t \le
T$.

Now set $U(x) = \sum_{1 \le i < j \le N}
\exp((\beta/3)\Phi(|x_i-x_j|))f(x_i)f(x_j)$, $x \in
\widetilde{(\Lambda_N)^N}$, $f \in C_c^1({\mathbb R}^d)$,
non-negative, and $\Phi$ as in condition (RP), (D). Then from
\eqref{diridomain} together with an approximation argument we can
conclude that $U \in D({\mathcal E}_{\scriptscriptstyle{\Lambda_N,
N}})$. This together with the fact that
$S^{\scriptscriptstyle{(\beta/3)}\Phi,
f}(\mbox{sym}^{\scriptscriptstyle{(N)}}(\cdot)) = U$ on
$\widetilde{(\Lambda_N)^N}$ implies via \eqref{dfconfi} that
$S^{\scriptscriptstyle{(\beta/3)}\Phi, f} \in D({\mathcal
E}^{\scriptscriptstyle{(N)}}_{\scriptscriptstyle{\Lambda_N}})$.
Hence as above we can find a $C_{10} < \infty$ such that for all
non-negative $f \in C^1_c({\mathbb R}^d), \, N \in {\mathbb N}$, $
\, 0 \le s \le t \le T$,
\begin{multline}\label{eq1606}
{\mathbb E}_{{\bf P}^{\scriptscriptstyle{(N)}}}
[|S^{\scriptscriptstyle{(\beta/3)}\Phi,f}({\bf  X}(t)) -
S^{\scriptscriptstyle{(\beta/3)}\Phi,f}({\bf X}(s))|^4]^{1/4} \le
C_{10} \, (t-s)^{1/2} \Big{(}
\int_{\Gamma_{\scriptscriptstyle{\Lambda_N}}^{\scriptscriptstyle{(N)}}}
G(\gamma) \,d\mu^{\scriptscriptstyle{(N)}}(\gamma) \Big{)}^{1/4},
\end{multline}
where
\begin{multline*}
G(\gamma) = \sum_{x \in \gamma} \sum_{y \in \gamma \setminus
\{x\}} \sum_{z \in \gamma \setminus \{x\}}
\exp((\beta/3)\Phi(|x-y|))
\exp((\beta/3)\Phi(|x-z|))  f(y) f(z) \\
\times \bigg(
\frac{\beta\Phi'(|x-y|)}{3|x-y|}\frac{\beta\Phi'(|x-z|)}{3|x-z|}
(x-y,x-z)_{{\mathbb R}^d} f(x)^2 + |\nabla f(x)|_{{\mathbb R}^d}^2  \\
+ \frac{\beta\Phi'(|x-y|)}{3|x-y|} (x-y, \nabla f(x))_{{\mathbb
R^d}} f(x) + \frac{\beta\Phi'(|x-z|)}{3|x-z|} (x-z, \nabla
f(x))_{{\mathbb R^d}} f(x) \bigg) \\ = \sum_{\{x,y\} \in \gamma}
g^f_2(x,y) + \sum_{\{x,y,z\} \in \gamma} g^f_3(x,y,z).
\end{multline*}
The function $g^f_2$ is given by
\begin{multline*}
g^f_2(x,y) = \exp((2/3)\beta\Phi(|x-y|)) \Bigg(\frac{2\beta^2}{9}
\Phi'(|x-y|)^2 f(x)^2 f(y)^2 + |\nabla f(x)|_{{\mathbb R}^d}^2
f(y)^2 \\ + |\nabla f(y)|_{{\mathbb R}^d}^2 f(x)^2 + \frac{2 \beta
\Phi'(|x-y|)}{3|x-y|} (x-y,
f(y)^2f(x) \nabla f(x) -f(x)^2 f(y) \nabla f(y))_{{\mathbb R^d}} \Bigg) \\
\end{multline*}
and $g^f_3$ is the symmetrization of
\begin{multline*}
6\exp((\beta/3)\Phi(|x-y|)) \exp((\beta/3)\Phi(|x-z|))  f(y) f(z)
\\ \times \bigg(
\frac{\beta\Phi'(|x-y|)}{3|x-y|}\frac{\beta\Phi'(|x-z|)}{3|x-z|}
(x-y,x-z)_{{\mathbb R}^d} f(x)^2 \\ + |\nabla f(x)|_{{\mathbb
R}^d}^2 + \frac{\beta\Phi'(|x-y|)}{3|x-y|} (x-y, \nabla
f(x))_{{\mathbb R^d}} f(x) + \frac{\beta\Phi'(|x-z|)}{3|x-z|}
(x-z, \nabla f(x))_{{\mathbb R^d}} f(x) \bigg).
\end{multline*}
Now by \eqref{corre} together with Theorem \ref{ruelle} we get
for all non-negative $f \in C^1_c({\mathbb R}^d), \, N \ge N_0$, $0 \le s \le t \le
T$, the following estimate:
\begin{eqnarray}\label{eq1607}
{\mathbb E}_{{\bf P}^{\scriptscriptstyle{(N)}}}
[|S^{\scriptscriptstyle{(\beta/3)}\Phi,f}({\bf  X}(t)) -
S^{\scriptscriptstyle{(\beta/3)}\Phi,f}({\bf X}(s))|^4]^{1/4}
\le  C_{10} \, (t-s)^{1/2} R(f),
\end{eqnarray}
where
\begin{multline}\label{eq1608}
R(f) := \Bigg{(} \frac{\zeta^3}{3!}\, \int_{({\mathbb R}^d)^3} |g^f_3(x_1, x_2, x_3)|
\exp\Big( - \frac{2}{3}\sum_{1
\le i < j \le 3} \beta \phi(x_i-x_j) \Big) \,dx^{\otimes 3} \\ +
\frac{\zeta^2}{2!}\, \int_{({\mathbb R}^d)^2} |g^f_2(x_1, x_2)|
\exp\Big( - \beta \phi(x_1-x_2) \Big) \,dx^{\otimes 2} \Bigg{)}^{1/4}.
\end{multline}
The integrals in \eqref{eq1608} are finite due to the
differentiability and integrability properties of $f$ and (RP),
(D) ($\Phi' \exp(-a\Phi)$ is by assumption a bounded function for
all $a > 0$). Then for all non-negative $f \in C^1_c({\mathbb R}^d), \, N \ge N_0$, $0 \le s \le t \le
T$,
\begin{eqnarray*}
{\mathbb E}_{{\bf P}^{\scriptscriptstyle{(N)}}}
[|S^{\scriptscriptstyle{(\beta/3)}\Phi,f}({\bf  X}(t)) -
S^{\scriptscriptstyle{(\beta/3)}\Phi,f}({\bf X}(s))|^4]^{1/4} \le
C_{10}  R(f) \, (t-s)^{1/2}.
\end{eqnarray*}
But then by \eqref{eq1606} there exists $C_{11} < \infty$ such that
\begin{eqnarray}\label{eqc11}
{\mathbb E}_{{\bf P}^{\scriptscriptstyle{(N)}}}
[|S^{\scriptscriptstyle{(\beta/3)}\Phi,f}({\bf  X}(t)) -
S^{\scriptscriptstyle{(\beta/3)}\Phi,f}({\bf X}(s))|^4]^{1/4} \le
C_{11} R(f) \, (t-s)^{1/2}
\end{eqnarray}
for all non-negative $f \in C^1_c({\mathbb R}^d), \, N \in {\mathbb N}$, $0 \le s \le t \le T$.
If we now assume that
\begin{eqnarray*}
q_k = \inf\{1, 1/ I(f_k)\} > 0 \quad \mbox{and} \quad p_k = \inf\{1, 1/ R(h_k) \} > 0,
\end{eqnarray*}
and $C_6 = C_9 + C_{11}$, then from \eqref{kolching} together with \eqref{eqc9} and
\eqref{eqc11} we can conclude \eqref{kolche}.
\hfill $\blacksquare$

{\bf Proof of Theorem \ref{tight}:} Criteria for tightness of
c\`adl\`ag (i.e., right continuous on $[0,\infty)$ and left limits
on $(0,\infty)$) processes in metric spaces have been worked out
in \cite[Chap.~3]{EK86}. For continuous processes as we are
considering one uses a slightly different modulus of continuity
(and also a different topology on the path space) as for
c\`adl\`ag processes. However, by using the Arzela--Ascoli Theorem
in metric spaces, see e.g.~\cite[Chap.~I, Theo.~23.2]{Choq66},
it is easy to show that \cite[Chap.~3, Theo.~7.2]{EK86} is also valid in the continuous
case. Since the sets
\begin{eqnarray*}
\{\gamma \in \Gamma \,|\, S_{\Phi,h}(\gamma) \le R \}, \qquad R <
\infty,
\end{eqnarray*}
are relatively compact subsets of $(\Gamma,
d_{\scriptscriptstyle{(\beta/3)} \Phi,h})$, see Proposition
\ref{propmetric}, Lemma \ref{tight0} yields condition (a) of
\cite[Chap.~3, Theo.~7.2]{EK86} (recall that
${\mu}^{\scriptscriptstyle{(N)}}$ is the invariant measure of
${\bf P}^{\scriptscriptstyle{(N)}}$). Condition (b) of
\cite[Chap.~3, Theo.~7.2]{EK86} follows from Lemma \ref{tight1}.
\hfill $\blacksquare$

\subsection{Identification of the limiting equilibrium
measures as a canonical Gibbs measures}

Consider the sequence of equilibrium measures
$(\mu^{\scriptscriptstyle{(N)}})_{N \in {\mathbb N}}$
corresponding to the $({\bf P}^{\scriptscriptstyle{(N)}})_{N \in
{\mathbb N}}$ as in Theorem \ref{tight}. Then tightness of $({\bf
P}^{\scriptscriptstyle{(N)}})_{N \in {\mathbb N}}$ implies
tightness of $(\mu^{\scriptscriptstyle{(N)}})_{N \in {\mathbb
N}}$. Now let $\mu$ be an accumulation point of
$(\mu^{\scriptscriptstyle{(N)}})_{N \in {\mathbb N}}$. Our aim is
to identify $\mu$ as a canonical Gibbs measure via an integration
by parts formula.

\begin{lemma}\label{conti}
{\rm Assume condition (D). For $n \in {\mathbb N}$ and
$v \in C^{\infty}_c({\mathbb R}^d, {\mathbb R}^d)$
consider the function
\begin{eqnarray*}
\Gamma \ni \gamma \mapsto L_{v,k}^{\phi}(\gamma) :=
- \beta \sum_{\{x,y\} \in \gamma}
(\nabla\phi(x-y),I_k(y) v(x)- I_k(x)v(y))_{{\mathbb R}^d},
\end{eqnarray*}
where the collection ${\bf I} = \{I_k \,|\, k \in {\mathbb N}\}$
is as in Section \ref{s1}. Then $L_{v,k}^{\phi}$ is a continuous
function on $(\Gamma, d_{\scriptscriptstyle{(\beta/3)}\Phi,h})$.}
\end{lemma}
{\bf Proof:} Just an easy modification of the proof of
\cite[Lem.~3.4]{KoKu04} where the continuity of
$S^{\scriptscriptstyle{(\beta/3)}\Phi,h_k}$ is shown. \hfill
$\blacksquare$

\begin{lemma}\label{l2limit}
{\rm Let the conditions in Theorem \ref{tight} hold. Then for all
accumulation points $\mu$ of $(\mu^{\scriptscriptstyle{(N)}})_{N
\in {\mathbb N}}$ and all $v \in C^{\infty}_c({\mathbb R}^d,
{\mathbb R}^d)$ we have that
\begin{eqnarray*}
L_{v}^{\phi, \mu} := \lim_{k \to \infty} L_{v,k}^{\phi}
\end{eqnarray*}
exists in $L^2(\mu)$.}
\end{lemma}
{\bf Proof:}
Set $f_k(x,y) := \beta \, \nabla\phi(x-y))\Big(I_k(y)v(x)-I_k(x)v(y)\Big)$,
$\{x, y\} \subset {\mathbb R}^d$. Then
\begin{eqnarray*}
(L_{v,k}^{\phi}(\gamma))^2 = \sum_{n=2}^4 \sum_{\{x_1, \ldots, x_n\} \subset \gamma} 
f^{[n,2]}_k(x_1, \ldots, x_n).
\end{eqnarray*}
Now as in the proof of Lemma \ref{tight0}, \eqref{corre} together with \eqref{expRuellebo}
yields for $N \ge N_0$ (as in Theorem \ref{ruelle})
\begin{eqnarray*}
{\mathbb E}_{{\mu}^{\scriptscriptstyle{(N)}}}
[(L_{v,k}^{\phi})^2]
\le  J_2 + J_3 + J_4,
\end{eqnarray*}
where
\begin{align*}
J_n = &\frac{\zeta^n}{n!} \int_{({\mathbb R}^d)^n} |f_k|^{[n,2]}(x_1, \ldots, x_n) 
\exp\Big( - \frac{2}{n}\sum_{1
\le i < j \le n} \beta \phi(x_i-x_j) \Big) \,dx^{\otimes n}, \quad n = 2, 3, 4.
\end{align*}
Since the potential $\phi$ is bounded from below, there exits
$C_{12} < \infty$ such that
\begin{multline*}
|J_4| \le  \frac{C_{12}}{4}
\Big(\int_{({\mathbb R}^d)^2} \|I_k(x_2) v(x_1)-I_k(x_1) v(x_2)\|_{{\mathbb R}^d}
\\ \times \|\beta \, \nabla\phi(x_1-x_2)\|_{{\mathbb R}^d}
\exp\Big( - \beta \phi(x_1-x_2) \Big) \,dx^{\otimes 2} \Big)^2
\le C_{12} \|v\|^2_{L^1(dx)} \|\beta \, \nabla \phi\|^2_{L^1(\exp(-\beta\phi)dx)}.
\end{multline*}
Analogously, (using Young's inequality) we obtain
\begin{multline*}
|J_3| \le \frac{C_{13}}{2} \|v\|_{\sup} \|\beta \, \nabla \phi\|_{L^1(\exp(-\beta\phi)dx)}
\\ \times \int_{({\mathbb R}^d)^2} \|I_k(x_2) v(x_1)-I_k(x_1) v(x_2)\|_{{\mathbb R}^d}
\|\beta \, \nabla\phi(x_1-x_2)\|_{{\mathbb R}^d}
\exp\Big( - \beta \phi(x_1-x_2) \Big) \,dx^{\otimes 2}
\\ \le C_{13} \|v\|_{\sup} \|v\|_{L^1(dx)} \|\beta \, \nabla \phi\|^2_{L^1(\exp(-\beta\phi)dx)}
\end{multline*}
for some $C_{13} < \infty$. $J_2$ can be estimated as follows:
\begin{multline*}
|J_2| \le \frac{C_{14}}{4} \int_{({\mathbb R}^d)^2} \|I_k(x_2) v(x_1)-I_k(x_1) v(x_2)\|^2_{{\mathbb R}^d}
\\ \times \|\beta \, \nabla\phi(x_1-x_2)\|^2_{{\mathbb R}^d}
\exp\Big( - \beta \phi(x_1-x_2) \Big) \,dx^{\otimes 2}
\le C_{14} \|v\|^2_{L^2(dx)} \|\beta \, \nabla \phi\|^2_{L^2(\exp(-\beta\phi)dx)}.
\end{multline*}
Next for $0 < r < \infty$ set
\begin{eqnarray*}
L_{v,k,r}^{\phi} := (L_{v,k,r}^{\phi} \vee -r) \wedge r.
\end{eqnarray*}
Then, by Lemma \ref{conti}, $L_{v,k,r}^{\phi}$ is a bounded
continuous function on $(\Gamma,
d_{\scriptscriptstyle{(\beta/3)}\Phi,h})$. Additionally,
\begin{eqnarray*}
(L_{v,k,r}^{\phi})^2 \le (L_{v,k}^{\phi})^2 \quad \mbox{and} \quad
(L_{v,k,r}^{\phi})^2 \nearrow (L_{v,k}^{\phi})^2  \quad \mbox{as} \quad
r \nearrow \infty.
\end{eqnarray*}
Now let $\mu$ be an accumulation point of
$(\mu^{\scriptscriptstyle{(N)}})_{N \in {\mathbb N}}$, i.e.,
$\mu^{\scriptscriptstyle{(N_n)}} \to \mu$ weakly for some subsequence
$N_n \to \infty$ as $n \to \infty$. Then
\begin{eqnarray*}
{\mathbb E}_{{\mu}}
[(L_{v,k,r}^{\phi})^2] = \lim_{n \to \infty}
{\mathbb E}_{{\mu}^{\scriptscriptstyle{(N_n)}}}[(L_{v,k,r}^{\phi})^2]
\le \liminf_{n \to \infty}
{\mathbb E}_{{\mu}^{\scriptscriptstyle{(N_n)}}}[(L_{v,k}^{\phi})^2]
\end{eqnarray*}
and so
\begin{eqnarray*}
{\mathbb E}_{{\mu}}
[(L_{v,k}^{\phi})^2] = \lim_{r \to \infty}
{\mathbb E}_{{\mu}}
[(L_{v,k,r}^{\phi})^2]
\le \liminf_{n \to \infty}
{\mathbb E}_{{\mu}^{\scriptscriptstyle{(N_n)}}}[(L_{v,k}^{\phi})^2] < \infty
\end{eqnarray*}
due to the estimates for $|J_2|, |J_3|, |J_4|$. Hence
$L_{v,k}^{\phi} \in L^2(\mu)$. Let $k \ge l$. As above we can
estimate
\begin{multline}\label{C16}
{\mathbb E}_{{\mu}}
[(L_{v,k}^{\phi} - L_{v,l}^{\phi})^2]
\le \liminf_{n \to \infty}
{\mathbb E}_{{\mu}^{\scriptscriptstyle{(N_n)}}}[(L_{v,k}^{\phi}
- L_{v,l}^{\phi})^2]
\le C_{15} \Bigg{(} \|v\|_{\sup} \|\beta \, \nabla \phi\|_{L^1(\exp(-\beta\phi)dx)}
\\ \times \int_{({\mathbb R}^d)^2} \|(I_k(x_2) - I_l(x_2)) v(x_1)-(I_k(x_1) - I_l(x_1))v(x_2)\|_{{\mathbb R}^d} \\
\times \|\beta \, \nabla\phi(x_1-x_2)\|_{{\mathbb R}^d}
\exp\Big( - \beta \phi(x_1-x_2) \Big) \,dx^{\otimes 2} \\
+\int_{({\mathbb R}^d)^2} \|(I_k(x_2) - I_l(x_2)) v(x_1)-(I_k(x_1) - I_l(x_1))v(x_2)\|^2_{{\mathbb R}^d} \\
\times \|\beta \, \nabla\phi(x_1-x_2)\|^2_{{\mathbb R}^d}
\exp\Big( - \beta \phi(x_1-x_2) \Big) \,dx^{\otimes 2} \Bigg{)} \\
\le C_{15} \Bigg{(}
2 \|v\|_{\sup} \|v\|_{L^1(dx)} \|\beta \, \nabla \phi\|_{L^1(\exp(-\beta\phi)dx)}
\\ \times \int_{{\mathbb R}^d \setminus B_{R_l}(0)}
 \|\beta \, \nabla\phi(x_1-x_2)\|_{{\mathbb R}^d}
\exp\Big( - \beta \phi(x_1-x_2) \Big) \,dx \\
+2 \|v\|^2_{L^2(dx)} \int_{{\mathbb R}^d \setminus B_{R_l}(0)} \|\beta \, \nabla\phi(x_1-x_2)\|^2_{{\mathbb R}^d}
\exp\Big( - \beta \phi(x_1-x_2) \Big) \,dx \Bigg{)} := C_{16}(l),
\end{multline}
where $R_l$ is a certain radius and $B_{R_l}(0)$ the corresponding
ball centered at the origin. Since $I_k(x) - I_l(x) = 0$ if
$\|x\| \le l-1$, we have $R_l \to \infty$ as $l \to \infty$. Now property (D) yields that
$C_{16}(l) \to 0$ as $l \to \infty$. Hence $(L_{v,k}^{\phi})_{k \in {\mathbb N}}$
is a Cauchy sequence in $L^2(\mu)$.
\hfill $\blacksquare$

\begin{lemma}\label{pointwiselimit}
{\rm Let the conditions in Theorem \ref{tight} hold, let $v \in
C^{\infty}_c({\mathbb R}^d, {\mathbb R}^d)$ and define
\begin{align*}
L_{v}^{\phi}(\gamma) := \left\{ \begin{array}{cc}
\!\!\!\!\!\!\!\!\!\!\!\!\!\!\!\!\!\!\!\!\!\!\!\!\!\!\!\!\!\!\!\!\!
\!\!\!\!\!\!\!\!\!\!\!-\beta \sum_{\{x,y\} \in \gamma}
(\nabla\phi(x-y),v(x)-v(y))_{{\mathbb R}^d} & \\
\quad\quad\quad\quad\quad\quad\quad\quad  \mbox{if} \sum_{\{x,y\}
\in \gamma} |(\nabla\phi(x-y), v(x))_{{\mathbb R}^d}| < \infty
\\\!\!\!\!\!\!\!\!\!\!\!\!\!\!\!\!\!\!\!\!\!\!\!\!\!\!\!\!\!\!\!\!\!
\!\!\!\!\!\!\!\!\!\!\!\!\!\!\!\!\!\!\!\!\!\!
\!\!\!\!\!\!\!\!\!\!\!\!\!\!\!\!\!\!\!\!\!\!\!\!\!\!\!\!\!\!\!\!\!
 0 & \!\!\!\!\!\!\!\!\!\!\!\!\!\!\!\!\!\!\!\!\!\!
 \!\!\!\!\!\!\!\!\!\!\!\!\!\!\!\!\!\!\!\!\!\!\!\!\!\!\!\!\!\!\!\!\!
\!\!\!\!\!\!\!\!\!\!\!\!\!\!\!\!\!\!\!\!\!\!\!\!\!\!\!\!\!\!\!\!\!
\!\!\!\!\!\!\!\!\!\!\!\!\!\!\!\!\!\!\!\!\!\!\mbox{otherwise}
\end{array}\right..
\end{align*}
Then $L_{v}^{\phi}$ is an $L^2(\mu)$-version of
$L_{v}^{\phi,\mu}$.}
\end{lemma}
{\bf Proof:} Define the sequence
\begin{eqnarray*}
M_{v,k}^{\phi}(\gamma) = \sum_{\{x,y\} \in \gamma}
|(\nabla\phi(x-y),v(x))_{{\mathbb R}^d}|I_k(y), \quad k \in
{\mathbb N}.
\end{eqnarray*}
Then $M_{v,k}^{\phi}(\gamma)$ monotonically converges to
\begin{eqnarray*}
M_{v}^{\phi}(\gamma) := \sum_{\{x,y\} \in \gamma}
|(\nabla\phi(x-y),v(x))_{{\mathbb R}^d}|
\end{eqnarray*}
as $k \to \infty$. Furthermore, by estimates as in the proof of
Lemma \ref{l2limit}, the $L^1(\mu)$-norms of the
$M_{v,k}^{\phi}(\gamma)$ are uniformly bounded. Thus, by monotone
convergence $M_{v}^{\phi} \in L^1(\mu)$ and therefore there exists
$S \subset \Gamma$ with $\mu(S) = 1$ such that
\begin{eqnarray*}
M_{v}^{\phi}(\gamma) < \infty \quad \mbox{for all} \quad \gamma
\in S.
\end{eqnarray*}

Now we return to the $L_{v,k}^{\phi}$. Note that for a subsequence
$L_{v,k_m}^{\phi}(\gamma) \to L_{v}^{\phi,\mu}(\gamma)$ as $m \to
\infty$ for $\mu$-a.a.~$\gamma \in \Gamma$. Obviously, for this
subsequence and  all $\gamma \in S$: $L_{v,k_m}^{\phi}(\gamma) \to
L_{v}^{\phi}(\gamma)$ as $m \to \infty$. Thus
\begin{eqnarray*}
L_{v}^{\phi,\mu}(\gamma) = \lim_{m \to \infty}
L_{v,k_m}^{\phi}(\gamma) = L_{v}^{\phi}(\gamma) \quad \mbox{for
$\mu$ a.a.} \quad \gamma \in S.
\end{eqnarray*}
\hfill $\blacksquare$

For later use we also need:

\begin{lemma}\label{l3norm}
{\rm Let the conditions in Theorem \ref{tight} hold. Then for all
subsequences $(\mu^{\scriptscriptstyle{(N_n)}})_{n \in {\mathbb N}}$ converging
weakly to an accumulation point $\mu$ of $(\mu^{\scriptscriptstyle{(N)}})_{N
\in {\mathbb N}}$ and all $v \in C^{\infty}_c({\mathbb R}^d,
{\mathbb R}^d)$ we have
\begin{eqnarray*}
\sup_{k \in {\mathbb N}} 
{\mathbb E}_{{\mu}}[|L_{v,k}^{\phi}|^3]
\le \sup_{k \in {\mathbb N}} \sup_{n \in {\mathbb N}} 
{\mathbb E}_{{\mu}^{\scriptscriptstyle{(N_n)}}}[|L_{v,k}^{\phi}|^3] < \infty.
\end{eqnarray*}
}
\end{lemma}
{\bf Proof:} As in the proof of Lemma \ref{l2limit} we get
for $N \ge N_0$ (as in Theorem \ref{ruelle}):
\begin{eqnarray*}
{\mathbb E}_{{\mu}^{\scriptscriptstyle{(N)}}}
[|L_{v,k}^{\phi}|^3]
\le  K_2 + K_3 + K_4 + K_5 + K_6,
\end{eqnarray*}
where
\begin{align*}
K_n = &\frac{\zeta^n}{n!} \int_{({\mathbb R}^d)^n} |f_k|^{[n,3]}(x_1, \ldots, x_n) 
\exp\Big( - \frac{2}{n}\sum_{1
\le i < j \le n} \beta \phi(x_i-x_j) \Big) \,dx^{\otimes n}, \quad n = 2, \ldots, 6.
\end{align*}
$K_2$ can be estimated as $J_2$, here we need that $\nabla \phi
\in L^3({\mathbb R}^d, \exp(-\beta\phi)dx)$.
$K_6$ can be estimated as $J_4$, here we need that $\nabla \phi
\in L^1({\mathbb R}^d, \exp(-\beta\phi)dx)$. Using Young's inequality $K3$-$K5$
can be treated as $J3$. In these cases we need $\nabla \phi
\in L^1({\mathbb R}^d, \exp(-\beta\phi)dx) \cap 
L^2({\mathbb R}^d, \exp(-\beta\phi)dx)$. Then as in the proof of Lemma \ref{l2limit} we get
the desired estimate.
\hfill $\blacksquare$

In order to formulate the next lemma we recall the gradient
$\nabla^\Gamma$ introduced and studied in \cite{AKR98a}. It acts
on finitely based smooth functions as in \eqref{fcbinfty} as
follows:
\begin{eqnarray*}
(\nabla^\Gamma F)(\gamma, x) = \sum_{j=1}^{n} \partial_j g_F
(\langle f_1, \gamma \rangle, \ldots, \langle f_n, \gamma \rangle)
\nabla f_j(x), \quad \gamma \in \Gamma, \, x \in \gamma.
\end{eqnarray*}
The corresponding directional derivative $\nabla_v^\Gamma$ in
direction $v \in C^{\infty}_c({\mathbb R}^d, {\mathbb R}^d)$ is
given by
\begin{eqnarray}\label{dderive}
(\nabla_v^\Gamma F)(\gamma) = \sum_{j=1}^{n} \partial_j g_F
(\langle f_1, \gamma \rangle, \ldots, \langle f_n, \gamma \rangle)
\langle (\nabla f_j, v)_{{\mathbb R}^d}, \gamma \rangle \quad
\gamma \in \Gamma, \, x \in \gamma.
\end{eqnarray}

\begin{lemma}\label{ipbfd}
{\rm Suppose that the conditions {\rm (BB), (D)} are satisfied and
that $N \in {\mathbb N}$, $\Lambda \subset {\mathbb R}^d$ bounded
Borel measurable. Let
$\mu^{\scriptscriptstyle{(N)}}_{\scriptscriptstyle{\Lambda}}$ be
the corresponding canonical Gibbs measure. Then for all $F, G \in
{\mathcal F} C_{\mathrm b}^\infty({\mathcal D},\Gamma)$ and $v \in
C^{\infty}_c(\mathring{\Lambda}, {\mathbb R}^d)$ the following
integration by parts formula holds:
\begin{eqnarray}\label{fdibp}
\int_{\Gamma^{\scriptscriptstyle{(N)}}_{\scriptscriptstyle{\Lambda}}}
\nabla_v^\Gamma F G
\,d\mu^{\scriptscriptstyle{(N)}}_{\scriptscriptstyle{\Lambda}}
=
- \int_{\Gamma^{\scriptscriptstyle{(N)}}_{\scriptscriptstyle{\Lambda}}}
F \nabla_v^\Gamma G
\,d\mu^{\scriptscriptstyle{(N)}}_{\scriptscriptstyle{\Lambda}}
- \int_{\Gamma^{\scriptscriptstyle{(N)}}_{\scriptscriptstyle{\Lambda}}}
F G B_v^{\phi}
\,d\mu^{\scriptscriptstyle{(N)}}_{\scriptscriptstyle{\Lambda}},
\end{eqnarray}
where
\begin{eqnarray}\label{bofdiv}
B_v^{\phi} := \langle \mbox{div} \,v, \cdot \rangle + L_v^{\phi}.
\end{eqnarray}}
\end{lemma}
{\bf Proof:} First let us show that $\langle \mbox{div} \,v, \cdot \rangle \in
L^2(d\mu^{\scriptscriptstyle{(N)}}_{\scriptscriptstyle{\Lambda}})$. Indeed, by \eqref{corre}
\begin{multline}\label{divv}
\int_{\Gamma^{\scriptscriptstyle{(N)}}_{\scriptscriptstyle{\Lambda}}}
\langle \mbox{div} \,v, \gamma \rangle^2
\,d\mu^{\scriptscriptstyle{(N)}}_{\scriptscriptstyle{\Lambda}}(\gamma)
= \int_{\Gamma^{\scriptscriptstyle{(N)}}_{\scriptscriptstyle{\Lambda}}}
\sum_{\{x,y\} \in \gamma} \mbox{div} \,v(x)\,\mbox{div} \,v(y)
+  \sum_{x \in \gamma} (\mbox{div} \,v(x))^2
\,d\mu^{\scriptscriptstyle{(N)}}_{\scriptscriptstyle{\Lambda}}(\gamma) \\
= \int_{\Lambda} (\mbox{div} \,v(x))^2
k_{\scriptscriptstyle{\Lambda}}^{\scriptscriptstyle{(1,N)}}(x) \,dx_\Lambda
+ \frac{1}{2}\, \int_{\Lambda^2}  \mbox{div} \,v(x_1)\,\mbox{div} \,v(x_2)
\, 
k_{\scriptscriptstyle{\Lambda}}^{\scriptscriptstyle{(2,N)}}
(x_1, x_2) \,dx_\Lambda^{\otimes 2}
\end{multline}
which is finite due to the boundedness of the correlation
functions. Similarly, we find that $\langle (\nabla f,
v)_{{\mathbb R}^d}, \cdot \rangle \in
L^2(d\mu^{\scriptscriptstyle{(N)}}_{\scriptscriptstyle{\Lambda}})$
for all $f \in {\mathcal D}$.

Next note that
\begin{eqnarray*}
B_v^{\phi}(\gamma) = \langle \mbox{div} \,v, \gamma \rangle \,\, +
L_{v,k}^{\phi}(\gamma)
\end{eqnarray*}
for all $\gamma \in
\Gamma^{\scriptscriptstyle{(N)}}_{\scriptscriptstyle{\Lambda}}$ if
$k$ is chosen large enough. Then using the ideas as in the proof
of Lemma \ref{l2limit} one gets $L_{v,k}^{\phi} \in
L^2(\mu^{\scriptscriptstyle{(N)}}_{\scriptscriptstyle{\Lambda}})$
and then, of course, also $B_{v}^{\phi} \in
L^2(\mu^{\scriptscriptstyle{(N)}}_{\scriptscriptstyle{\Lambda}})$.
Now by going to Euclidean coordinates one easily proves
\eqref{fdibp} by integrating by parts. Note that the boundary
terms are zero due to the support property of $v$ and that (D)
implies that
\begin{eqnarray*}
\nabla \exp(-\phi) = - \nabla \phi \exp(-\phi) \quad dx\mbox{-a.e.~on} \,\, {\mathbb R}^d.
\end{eqnarray*}
\hfill $\blacksquare$

\begin{theorem}\label{thidibp}
Assume the conditions in Theorem \ref{tight}. Furthermore, let
$\mu$ be an accumulation point of
$(\mu^{\scriptscriptstyle{(N)}})_{N \in {\mathbb N}}$ provided by
Theorem \ref{tight}. Then for all $F, G \in {\mathcal F}
C_{\mathrm b}^\infty({\mathcal D},\Gamma)$ and $v \in
C^{\infty}_c({\mathbb R}^d, {\mathbb R}^d)$ the following
integration by parts formula holds:
\begin{eqnarray}\label{idibp}
\int_{\Gamma} \nabla_v^\Gamma F G \,d\mu = - \int_{\Gamma} F
\nabla_v^\Gamma G \,d\mu - \int_{\Gamma} F G B_v^{\phi} \,d\mu.
\end{eqnarray}
In particular, $\mu$ is a canonical Gibbs measure.
\end{theorem}
{\bf Proof:} By the product rule for $\nabla^\Gamma$ it suffices
to prove \eqref{idibp} for
\begin{eqnarray*}
F =g_F(\langle f_1,\cdot\rangle,\dots,\langle f_n,\cdot\rangle)
\quad \mbox{and} \quad G \equiv 1.
\end{eqnarray*}
If now $\mu^{\scriptscriptstyle{(N_n)}} \to \mu$ weakly as $n \to
\infty$, then by Lemma \ref{ipbfd} it suffices to show that
\begin{eqnarray}\label{bigtask}
\lim_{n \to \infty} {\mathbb E}_{\mu^{\scriptscriptstyle{(N_n)}}}
[\nabla_v^\Gamma F] = {\mathbb E}_{\mu}[\nabla_v^\Gamma F] \quad
\mbox{and} \quad \lim_{n \to \infty} {\mathbb
E}_{\mu^{\scriptscriptstyle{(N_n)}}} [F B_v^{\phi}] = {\mathbb
E}_{\mu}[F B_v^{\phi}].
\end{eqnarray}

Let us first consider the second identity in \eqref{bigtask}. From
\eqref{divv} together with \eqref{Ruellebo} we can conclude that
\begin{eqnarray}\label{C17}
{\mathbb E}_{{\mu}^{\scriptscriptstyle{(N)}}}[\langle \mbox{div} \,v, \cdot \rangle^2]
\le \xi \|\mbox{div} \,v \|^2_{L^2(dx)} +
\frac{\xi^2}{2} \|\mbox{div} \,v \|^2_{L^1(dx)} := C_{17}.
\end{eqnarray}
Furthermore, notice that for $0 < r < \infty$
\begin{eqnarray*}
D_{v,r} := (\langle \mbox{div} \,v, \cdot \rangle \vee -r) \wedge r
\end{eqnarray*}
is a bounded continuous function on $(\Gamma,
d_{\scriptscriptstyle{(\beta/3)}\Phi,h})$. Hence
\begin{eqnarray*}
{\mathbb E}_{\mu}[D_{v,r}^2]
\le {\mathbb E}_{\mu}[\langle \mbox{div} \,v, \cdot \rangle^2]
\le \xi \|\mbox{div} \,v \|^2_{L^2(dx)} +
\frac{\xi^2}{2} \|\mbox{div} \,v \|^2_{L^1(dx)},
\end{eqnarray*}
by the same arguments as in the proof of Lemma \ref{l2limit}
(there applied to $L_{v,k}^{\phi}$). Then by the triangle
inequality, \eqref{C16} and \eqref{C17}
\begin{multline*}
|{\mathbb E}_{\mu^{\scriptscriptstyle{(N_n)}}} [F B_v^{\phi}]
-{\mathbb E}_{\mu}[F B_v^{\phi}]| \le {\mathbb
E}_{\mu^{\scriptscriptstyle{(N_n)}}}
[|F||\langle \mbox{div} \,v, \cdot \rangle - D_{v,r}|] \\
+|{\mathbb E}_{\mu^{\scriptscriptstyle{(N_n)}}} [FD_{v,r})] -
{\mathbb E}_{\mu}[FD_{v,r})]| + {\mathbb E}_{\mu}
[|F||\langle \mbox{div} \,v, \cdot \rangle - D_{v,r}|] \\
+ {\mathbb E}_{\mu^{\scriptscriptstyle{(N_n)}}} \Big[|F|\Big|
\sum_{\{x,y\} \in (\cdot)}
(\nabla\phi(x-y), (1-I_k(y))v(x)-(1-I_k(x))v(y))_{{\mathbb R}^d}\Big| \Big] \\
+ {\mathbb E}_{\mu^{\scriptscriptstyle{(N_n)}}}
[|F||L_{v,k}^{\phi} - L_{v,k,r}^{\phi}|] + |{\mathbb
E}_{\mu^{\scriptscriptstyle{(N_n)}}}
[F L_{v,k,r}^{\phi}] - {\mathbb E}_{\mu}[F L_{v,k,r}^{\phi}]| \\
+ {\mathbb E}_{\mu}[|F||L_{v,k,r}^{\phi} - L_{v,k}^{\phi}|]
+ {\mathbb E}_{\mu}[|F||L_{v,k}^{\phi} - L_{v}^{\phi}|] \\
\le \frac{2 C_{17} \|g_F\|_{\sup}}{r} +  |{\mathbb
E}_{\mu^{\scriptscriptstyle{(N_n)}}}
[FD_{v,r})] - {\mathbb E}_{\mu}[FD_{v,r})]| \\
+ \|g_F\|_{\sup} C_{16}(k) + |{\mathbb
E}_{\mu^{\scriptscriptstyle{(N_n)}}} [F L_{v,k,r}^{\phi}] -
{\mathbb E}_{\mu}[F L_{v,k,r}^{\phi}]|
+ {\mathbb E}_{\mu}[|F||L_{v,k}^{\phi} - L_{v}^{\phi}|] \\
+  \|g_F\|_{\sup}\frac{ \sup_{k \in {\mathbb N}}
{\mathbb E}_{{\mu}}[(L_{v,k}^{\phi})^2]}{r}
+ \|g_F\|_{\sup}\frac{ \sup_{k \in {\mathbb N}}\sup_{n \in {\mathbb N}}
{\mathbb E}_{{\mu}^{\scriptscriptstyle{(N_n)}}}[(L_{v,k}^{\phi})^2]}{r}.
\end{multline*}
The constants  $\sup_{k \in {\mathbb N}}{\mathbb
E}_{{\mu}}[(L_{v,k}^{\phi})^2]$,
$\sup_{k \in {\mathbb N}}\sup_{n \in {\mathbb N}} {\mathbb
E}_{{\mu}^{\scriptscriptstyle{(N_n)}}}[(L_{v,k}^{\phi})^2]$ are
finite due to the estimates for $|J_2|, |J_3|,$ $|J_4|$ in the
proof of Lemma \ref{l2limit}. Now the second identity in
\eqref{bigtask} follows from Lemma \ref{l2limit} and the weak
convergence ${\mu}^{\scriptscriptstyle{(N_n)}} \to \mu$ as $n \to
\infty$.

Note that
\begin{eqnarray*}
\nabla_v^\Gamma F = \sum_{j=1}^{n}
\partial_j g_F (\langle f_1, \cdot \rangle, \ldots, \langle f_n,
\cdot \rangle) \langle (\nabla f_j, v)_{{\mathbb R}^d}, \cdot
\rangle.
\end{eqnarray*}
Thus, showing the first identity in \eqref{bigtask} is a special
case of proving the second one. Just take the monomial $\langle
(\nabla f_j, v)_{{\mathbb R}^d}, \cdot \rangle$ instead of the
monomial $\langle \mbox{div} \,v, \cdot \rangle$ and the function
$1$ instead of $L_{v}^{\phi}$, all the other functions involved
are bounded and continuous.

Hence we have shown \eqref{idibp}. The fact that $\mu$ is a
canonical Gibbs measure now follows from \cite[Theo.~4.3]{AKR98b}.
\hfill $\blacksquare$

\subsection{Identification of the accumulation points
as the distribution of an infinite volume,
infinite particle stochastic dynamics}

Let us fix an accumulation point $\mu$ of
$(\mu^{\scriptscriptstyle{(N)}})_{N \in {\mathbb N}}$. Then for
all $F, G \in {\mathcal F}C_b^{\infty}({\mathcal D}, \Gamma)$ we
consider the bilinear from
\begin{multline}\label{dfongamma}
{\mathcal E}_{\mu}(F,G) : =
\int_{\Gamma} \langle \nabla^{\Gamma}F(\gamma),
\nabla^{\Gamma}G(\gamma) \rangle_{T_\gamma(\Gamma)} \,d \mu(\gamma) \\ 
= \int_{\Gamma} \sum_{x \in \gamma} (\nabla^{\Gamma}F(\gamma,x),
\nabla^{\Gamma}G(\gamma,x))_{{\mathbb R}^d} \,d \mu(\gamma),
\end{multline}
where $\langle \cdot, \cdot \rangle_{T_\gamma(\Gamma)} = \sum_{x \in \gamma} 
(\cdot, \cdot)_{{\mathbb R}^d}$ is the scalar product in the tangent space 
$T_\gamma(\Gamma)$, see \cite{AKR98b} for details. 
Using the integration by parts formula derived in Theorem \ref{thidibp} we obtain
for $F, G \in {\mathcal F}C_b^{\infty}({\mathcal D}, \Gamma)$:
\begin{eqnarray*}
{\mathcal E}_{\mu}(F,G) =
\int_{\Gamma} H_{\mu} F G \,d \mu,
\end{eqnarray*}
where
\begin{multline}\label{hmu}
H_{\mu}F = - \sum_{i,j=1}^{N} \partial_i
\partial_j g_F(\langle f_1, \cdot \rangle, \ldots, \langle f_N,
\cdot \rangle) \langle (\nabla f_i, \nabla f_j)_{{\mathbb R}^d},
\cdot \rangle \\
- \sum_{j=1}^{N} \partial_j g_F(\langle f_1, \cdot \rangle,
\ldots, \langle f_N, \cdot \rangle) \Big{(} \langle \Delta f_j,
\cdot \rangle  + L^{\phi}_{\nabla f_j}\Big{)}
\end{multline}
for $F\in {\mathcal F}C_{b}^\infty ({\mathcal D}, \Gamma )$ as in
\eqref{fcbinfty}.


\begin{theorem}\label{martingaleproblem}
Assume the conditions as in Theorem \ref{tight}. Furthermore, let
${\bf P}$ be an accumulation point of $({\bf
P}^{\scriptscriptstyle{(N)}})_{N \in {\mathbb N}}$ with invariant
canonical Gibbs measure $\mu$ provided in Theorem \ref{tight}.
Then ${\bf P}$ solves the martingale problem for $(-H_{\mu},
{\mathcal F}C_b^{\infty}({\mathcal D}, \Gamma))$ with initial
distribution $\mu$, i.e., for all $G \in {\mathcal
F}C_b^{\infty}({\mathcal D}, \Gamma)$,
\begin{eqnarray}\label{martprob}
G({\bf X}(t)) - G(({\bf X}(0)) + \int_0^t H_{\mu}G(({\bf X}(u))
\,du, \qquad t \ge 0,
\end{eqnarray}
is an ${\bf F}_t$-martingale under ${\bf P}$ and ${\bf P} \circ
{\bf X}(0)^{-1} = \mu$.
\end{theorem}
{\bf Proof:} For $t, s \ge 0$, we define the following random
variable on $C([0, \infty), \Gamma)$:
\begin{eqnarray*}
U({\bf X},t,s) := G({\bf X}(t+s)) - G({\bf X}(t)) + \int_t^{t+s}
H_{\mu} G({\bf X}(u))\,du.
\end{eqnarray*}
Corresponding to
\begin{eqnarray*}
G=g_G(\langle f_1,\cdot\rangle,\dots,\langle f_n,\cdot\rangle) \in
{\mathcal F}C_b^{\infty}({\mathcal D}, \Gamma)
\end{eqnarray*}
we define
\begin{eqnarray*}
\widetilde{G}(x) := g_G\Big(\sum_{i=1}^N f_1(x_i),\dots,
\sum_{i=1}^N f_n(x_i)\Big), \quad x = (x_1, \ldots, x_N) \in
\Lambda^N.
\end{eqnarray*}
Note that $G(\mbox{sym}^{\scriptscriptstyle{(N)}}(\cdot)) =
\widetilde{G}$ on $\widetilde{\Lambda^N}$. Since the $f_1, \ldots,
f_n$ have compact support there exists $N_0 \in {\mathbb N}$ such
that $\widetilde{G}$ is an element of
$C^2_c(\mathring{(\Lambda_N)^N}) \subset
D(H_{\scriptscriptstyle{\Lambda_N\!, N}})$ for all $N \ge N_0$.
Hence for $N \ge N_0$, $G \in
D(H^{\scriptscriptstyle{(N)}}_{\scriptscriptstyle{\Lambda_N}})$
and we have the pointwise representation of
$H^{\scriptscriptstyle{(N)}}_{\scriptscriptstyle{\Lambda_N}}G$
provided in \eqref{explicitgenconfi}. Notice that this
representation coincides with the pointwise representation of
$H_\mu G$, see \eqref{hmu}.

The trace filtration obtained by restricting $({\bf F}_t)_{t \ge
0}$ to $C([0, \infty),
\Gamma_{\scriptscriptstyle{\Lambda}}^{\scriptscriptstyle{(N)}})$
coincides with the natural filtration of $C([0, \infty),
\Gamma_{\scriptscriptstyle{\Lambda}}^{\scriptscriptstyle{(N)}})$.
Furthermore, ${\bf P}^{\scriptscriptstyle{(N)}}$ solves the
martingale problem for
$(-H^{\scriptscriptstyle{(N)}}_{\scriptscriptstyle{\Lambda_N}},
D(H^{\scriptscriptstyle{(N)}}_{\scriptscriptstyle{\Lambda_N}}))$
w.r.t.~$({\bf
F}^{\scriptscriptstyle{(N)}}_{\scriptscriptstyle{\Lambda_N}}(t))_{t
\ge 0}$. Therefore we have for all ${\bf F}_t$-measurable,
bounded, continuous $F_t: C([0, \infty), \Gamma)$ $\to {\mathbb
R}$ and $N \ge N_0$ that ${\mathbb E}_{{\bf
P}^{\scriptscriptstyle{(N)}}} [F_t U(t, s)] = 0$. Thus it follows
that
\begin{eqnarray}\label{expzero}
0 = \lim_{N \to \infty} {\mathbb E}_{{\bf
P}^{\scriptscriptstyle{(N)}}} [F_t U(t, s)].
\end{eqnarray}
Now let $(N_n)_{n \in {\mathbb N}}$ a subsequence such that ${\bf
P}^{\scriptscriptstyle{(N_n)}} \to {\bf P}$ weakly. Having
\eqref{expzero} it remains to show
\begin{eqnarray}\label{littletask}
\lim_{n \to \infty}{\mathbb E}_{{\bf
P}^{\scriptscriptstyle{(N_n)}}}[F_t U(t,s)] = {\mathbb E}_{{\bf
P}}[F_t U(t,s)]
\end{eqnarray}
to obtain \eqref{martprob}. We have
\begin{multline}\label{moredetail}
|{\mathbb E}_{{\bf P}^{\scriptscriptstyle{(N_n)}}}[F_t U(t,s)] -
{\mathbb E}_{{\bf P}}[F_t U(t,s)]| \\ \le |{\mathbb E}_{{\bf
P}^{\scriptscriptstyle{(N_n)}}}[F_t(G({\bf X}(t+s)) - G({\bf
X}(t)))] -
{\mathbb E}_{{\bf P}}[F_t(G({\bf X}(t+s)) - G({\bf X}(t)))]| \\
+ \int_t^{t+s} |{\mathbb E}_{{\bf
P}^{\scriptscriptstyle{(N_n)}}}[F_t H_{\mu} G({\bf X}(u))] -
{\mathbb E}_{{\bf P}}[F_t H_{\mu} G({\bf X}(u))]| \,du.
\end{multline}
The first term on the right hand side of the estimate
\eqref{moredetail} converges to zero as $n \to \infty$, because
the function $F_t(G({\bf X}(t+s)) - G({\bf X}(t)))$ is bounded and
continuous. Showing that
\begin{eqnarray*}
|{\mathbb E}_{{\bf P}^{\scriptscriptstyle{(N_n)}}}[F_t H_{\mu}
G({\bf X}(u))] - {\mathbb E}_{{\bf P}}[F_t H_{\mu} G({\bf X}(u))]|
\to 0 \quad \mbox{as} \quad n \to \infty \quad \forall u \in
[t,t+s]
\end{eqnarray*}
is essentially the same as proving \eqref{bigtask}, done in the
proof of Theorem \ref{thidibp}. Now using the Cauchy--Schwartz
inequality, the fact that $\mu$ and
$\mu^{\scriptscriptstyle{(N_n)}}$ are the invariant measures of
${\bf P}$ and ${\bf P}^{\scriptscriptstyle{(N_n)}}$, respectively,
and the boundedness of $\{ {\mathbb
E}_{{\mu}^{\scriptscriptstyle{(N_n)}}}[(H_{\mu} G)^2] \,|\, n \in
{\mathbb N}\}$ we find a constant $C_{18} < \infty$ independent of
$u \in [t,t+s]$ and $n \in {\mathbb N}$ such that
\begin{eqnarray*}
|{\mathbb E}_{{\bf P}^{\scriptscriptstyle{(N_n)}}}[F_t H_{\mu}
G({\bf X}(u))] - {\mathbb E}_{{\bf P}}[F_t H_{\mu} G({\bf X}(u))]|
\le C_{18}.
\end{eqnarray*}
Therefore, the second term on the right hand side of the estimate
\eqref{moredetail} converges to zero as $n \to \infty$ by Lebesgue
dominated convergence. Thus, \eqref{littletask} is shown.

Obviously, we have ${\bf P} \circ {\bf X}(t)^{-1} = \mu$ for all
$t \ge 0$, in particular ${\bf P} \circ {\bf X}(0)^{-1} = \mu$.
\hfill $\blacksquare$

\begin{remark}\label{solution}
{\rm Using It\^o's formula, Theorem \ref{martingaleproblem}
implies that each accumulation point ${\bf P}$ of $({\bf
P}^{\scriptscriptstyle{(N)}})_{N \in {\mathbb N}}$ solves the
following infinite system of stochastic differential equation in
the sense of the associated martingale problem:
\begin{align}\label{ivsd}
dx(t) & = -\beta \!\! \sum_{\stackunder{y(t) \neq x(t)} {y(t) \in
{\bf X}(t)}} \nabla \phi(x(t) - y(t)) \,dt + \sqrt{2}\,dB^x(t), &
 \nonumber \\
{\bf P} \circ {\bf X}(0)^{-1} & = \mu,
\end{align}
where $x(t) \in {\bf X}(t) \in \Gamma$, $(B^{x})_{x \in \gamma}$,
$\gamma \in \Gamma$, is a sequence of independent Brownian motions and $\mu$ is the
invariant measure corresponding to ${\bf P}$.}
\end{remark}

\subsection{Identification of the accumulation points
as Markov processes and uniqueness}

By \cite{AKR98b} the closure of $({\mathcal E}_\mu,
{\mathcal F}C_b^{\infty}({\mathcal D}, \Gamma))$ on $L^2(\Gamma,
\mu)$, in sequel denoted by $({\mathcal E}^{\rm min}_\mu,
D({\mathcal E}^{\rm min}_\mu))$, is conservative, local and quasi-regular, 
hence associated with a diffusion process on $\Gamma$. 
When started with $\mu$ its distribution ${\bf P}_\mu$ also satisfies 
the martingale problem \eqref{martprob}. So far we do not know 
whether ${\bf P}_\mu = {\bf P}$ with ${\bf P}$ as in
Theorem \ref{martingaleproblem}. A first step to that identification
yields the following convergence of the associated Dirichlet forms. 

\begin{proposition}\label{thmconvdf}
Let the assumptions in Theorem \ref{tight} hold and 
let $(\mu^{\scriptscriptstyle{(N_n)}
})_{n \in {\mathbb N}}$ be a subsequence
converging to an accumulation point $\mu$ of
$(\mu^{\scriptscriptstyle{(N)}})_{N \in {\mathbb N}}$. 
Then for all $F, G \in {\mathcal F}C_{\mathrm b}^\infty({\mathcal D},\Gamma)$ 
\begin{eqnarray}\label{convdf}
\lim_{n \to \infty} {\mathcal E}^{\scriptscriptstyle{(N_n)}}_{\scriptscriptstyle{\Lambda_{N_n}}}(F,G)
= {\mathcal E}_\mu(F,G).
\end{eqnarray}
\end{proposition}
{\bf Proof:} By polarization identity we can restrict ourself to 
the case $F=G$. From \eqref{dfongamma} we get that 
\begin{multline*}
{\mathcal E}_\mu(F,F) \\ = \sum_{i,j=1}^{n} 
\int_{\Gamma} \partial_i g_F
(\langle f_1, \gamma \rangle, \ldots, \langle f_n, \gamma \rangle)
\partial_j g_F
(\langle f_1, \gamma \rangle, \ldots, \langle f_n, \gamma \rangle)
\langle (\nabla f_i,\nabla f_j)_{{\mathbb R}^d}, \gamma \rangle \,d \mu(\gamma) 
\end{multline*}
for $F$ as in \eqref{fcbinfty}. Furthermore, by definition of 
${\mathcal E}^{\scriptscriptstyle{(N)}}_{\scriptscriptstyle{\Lambda}}$, see
\eqref{dfconfi}, we find
\begin{multline*}
{\mathcal E}^{\scriptscriptstyle{(N_n)}}_{\scriptscriptstyle{\Lambda_{N_n}}}(F,F) 
\\ = \sum_{i,j=1}^{n} 
\int_{\Gamma} \partial_i g_F
(\langle f_1, \gamma \rangle, \ldots, \langle f_n, \gamma \rangle)
\partial_j g_F
(\langle f_1, \gamma \rangle, \ldots, \langle f_n, \gamma \rangle)
\langle (\nabla f_i,\nabla f_j)_{{\mathbb R}^d}, \gamma \rangle 
\,d\mu^{\scriptscriptstyle{(N_n)}}(\gamma) 
\end{multline*}
again for $F$ as in \eqref{fcbinfty}.
Since $\mu^{\scriptscriptstyle{(N_n)}} \to \mu$ weakly as $n \to
\infty$, \eqref{convdf} follows by analogous arguments as in the proof of 
Theorem \ref{thidibp}.
\hfill$\blacksquare$

This convergence, however, is too weak to conclude convergence of the associated
semi-groups or resolvents. For this we need the stronger Mosco convergence 
of quadratic forms. The concepts of Mosco convergence were introduced in 
\cite{Mos94}. Here we need a generalization of these concepts provided in 
\cite{KuSh03}.     

\begin{definition}\label{defconhs}
{\rm We say that a sequence of Hilbert spaces $(H_n)_{n \in {\mathbb N}}$ converges
to a Hilbert space $H$, if there exists a dense subspace $C \subset H$ and a 
sequence of operators $(\Phi_n)_{n \in {\mathbb N}}$, where
\begin{eqnarray*}
\Phi_n: C \to H_n, \quad n \in {\mathbb N},
\end{eqnarray*}
with the following property:
\begin{eqnarray*}
\lim_{n \to \infty} \|\Phi_n u \|_{H_n} =  \|u \|_{H}
\end{eqnarray*}
for all $u \in C$.
}
\end{definition}

Let $\mu$ be an accumulation point of $(\mu^{\scriptscriptstyle{(N)}})_{N \in {\mathbb N}}$ 
and $(\mu^{\scriptscriptstyle{(N_n)}})_{n \in {\mathbb N}}$ a subsequence
such that  $\lim_{n \to \infty } \mu^{\scriptscriptstyle{(N_n)}} = \mu$.
When choosing $C := {\mathcal F}C_{\mathrm b}^\infty({\mathcal D},\Gamma)$ and
the mapping $\Phi_n := R_n$, $n \in {\mathbb N}$, as the choice of the continuous representative
of a function from ${\mathcal F}C_{\mathrm b}^\infty({\mathcal D},\Gamma) \subset
L^2(\mu)$ (this can be done uniquely, since $\mu$ as a Gibbs measure has full topological 
support on $\Gamma$) and then considering as function in $L^2(\mu^{\scriptscriptstyle{(N_n)}})$,
we see that $H_n := L^2(\mu^{\scriptscriptstyle{(N_n)}})$ converges to 
$H := L^2(\mu)$ in the sense of Definition \ref{defconhs} as $n \to \infty$.

\begin{definition}[strong convergence]\label{strong}
{\rm Let $(H_n)_{n \in {\mathbb N}}$, $(\Phi_n)_{n \in {\mathbb N}}$, $H$ 
and $C$ be as in Definition \ref{defconhs}. We say that a sequence of vectors 
$(u_n)_{n \in {\mathbb N}}$ with $u_n \in H_n, n \in {\mathbb N},$ converges 
strongly to a vector $u \in H$, if there exists a sequence  
$(\tilde{u}_n)_{n \in {\mathbb N}}$ in $C$ with the following properties:
\begin{align*}
\lim_{m \to \infty} &\|\tilde{u}_m - u\|_{H} = 0 \\
\lim_{m \to \infty} \limsup_{n \to \infty} &\|\Phi_n \tilde{u_m} - u_n\|_{H_n} = 0.
\end{align*}
}
\end{definition}


\begin{definition}[weak convergence]\label{weak}
{\rm Let $(H_n)_{n \in {\mathbb N}}$, $(\Phi_n)_{n \in {\mathbb N}}$, $H$ 
and $C$ be as in Definition \ref{defconhs}. We say that a sequence of vectors 
$(u_n)_{n \in {\mathbb N}}$ with $u_n \in H_n, n \in {\mathbb N},$ converges 
weakly to a vector $u \in H$, if 
\begin{align*}
\lim_{n \to \infty} (u_n, v_n)_{H_n} =  (u, v)_{H} 
\end{align*}
for every sequence $(v_n)_{n \in {\mathbb N}}$ 
with $v_n \in H_n$, $n \in {\mathbb N}$, which strongly converges to $v \in H$.
}
\end{definition}

In \cite{Kol05}[Lem.~2.7] the following simple criterion for strong convergence
has been proved.

\begin{lemma}\label{lemstrong}
{\rm Let $(H_n)_{n \in {\mathbb N}}$, $(\Phi_n)_{n \in {\mathbb N}}$, $H$ 
and $C$ be as in Definition \ref{defconhs}. A sequence 
$(u_n)_{n \in {\mathbb N}}$ with $u_n \in H_n, n \in {\mathbb N},$ converges 
strongly to a vector $u \in H$, if and only if
\begin{align*}
\lim_{n \to \infty} \|u_n\|_{H_n} =  \|u\|_{H} \quad \mbox{and} 
\quad \lim_{n \to \infty} (u_n, \Phi_n(v))_{H_n} =  (u, v)_{H} \quad
\mbox{for all} \quad v \in C.
\end{align*}
}
\end{lemma}

\begin{definition}\label{convoper}
{\rm Let $(H_n)_{n \in {\mathbb N}}$, $(\Phi_n)_{n \in {\mathbb N}}$, $H$ 
and $C$ be as in Definition \ref{defconhs}. We say that a sequence of bounded
operators $(B_n)_{n \in {\mathbb N}}$ with $B_n \in L(H_n), n \in {\mathbb N},$ converges 
strongly to a bounded operator $B \in L(H)$, if for every sequence 
$(u_n)_{n \in {\mathbb N}}$ with $u_n \in H_n$ , $n \in {\mathbb N}$, which strongly 
converges to $u \in H$, the sequence $(B_n u_n)_{n \in {\mathbb N}}$ strongly
converges to $Bu$.
}
\end{definition}

Next we consider convergence of quadratic forms ${\mathcal Q}$. Recall that a quadratic form 
on a Hilbert space ${\mathcal H}$ is given by a bilinear form 
${\mathcal E}: D({\mathcal E}) \times D({\mathcal E}) \to {\mathbb R}$, 
where $D({\mathcal E}) \subset {\mathcal H}$.  
We consider only densely defined, non-negative, closed, symmetric bilinear forms. Then we
define the corresponding quadratic form  ${\mathcal Q}: {\mathcal H} \to \overline{\mathbb R}$     
by setting  
\begin{align*}
{\mathcal Q}(u) := \left\{ \begin{array}{cc} {\mathcal E}(u,u) \quad \mbox{if} 
\quad u \in D({\mathcal E}) \\
\!\!\!\!\!\!\!\!\!\!\!\!\!\!\!\!\!\!\!\!\!\!\!\!\!\!\!\!\!\!
\!\!\!\!\!\!\infty &  \!\!\!\!\! \!\!\!\!\!
 \!\!\!\!\! \!\!\!\!\! \!\!\!\!\! \!\!\!\!\!
 \!\!\!\!\! \mbox{otherwise}
\end{array}\right..
\end{align*}
Recall that closedness of $({\mathcal E}, D({\mathcal E}))$ is equivalent to 
lower semi-continuity of ${\mathcal Q}: {\mathcal H} \to \overline{\mathbb R}$. 

\begin{definition}[Mosco convergence]\label{moscoconv}
{\rm Let $(H_n)_{n \in {\mathbb N}}$, $(\Phi_n)_{n \in {\mathbb N}}$, $H$ 
and $C$ be as in Definition \ref{defconhs}. We say that a sequence of quadratic
forms $({\mathcal Q}_n)_{n \in {\mathbb N}}$ with ${\mathcal Q}_n: H_n \to 
\overline{\mathbb R}$, $n \in {\mathbb N}$, Mosco converges 
to a quadratic form  ${\mathcal Q}: H \to 
\overline{\mathbb R}$, if the following conditions hold:
\begin{enumerate}
\item[(M1)] If a sequence $(u_n)_{n \in {\mathbb N}}$ with 
$u_n \in H_n, n \in {\mathbb N},$ weakly converges to a vector $u \in H$, then
\begin{eqnarray*}
{\mathcal Q} (u) \le \liminf_{n \to \infty} {\mathcal Q}_n(u_n).
\end{eqnarray*}

\item[(M2)] For all $u \in H$ there exists a sequence $(u_n)_{n \in {\mathbb N}}$
with $u_n \in H_n, n \in {\mathbb N},$ which strongly converges to $u$ and   
\begin{eqnarray*}
{\mathcal Q} (u) = \lim_{n \to \infty} {\mathcal Q}_n(u_n).
\end{eqnarray*}
\end{enumerate}
}
\end{definition}

In \cite{Mos94} it is proved that Mosco convergence of a sequence of quadratic
forms is equivalent to the convergence, in the strong operator sense, of the 
sequence of semi-groups and resolvents, respectively, associated with the 
corresponding bilinear forms. In \cite{KuSh03} this result is generalized to 
the present situation, where we have a sequence of Hilbert spaces. Here strong 
convergence of bounded operators has to be understood in the sense of Definition 
\ref{convoper}.      
 
We are interested in the case where ${\mathcal Q}_n$ is the quadratic form corresponding to
$({\mathcal E}^{\scriptscriptstyle{(N_n)}}_{\scriptscriptstyle{\Lambda_{N_n}}}$, 
$D({\mathcal E}^{\scriptscriptstyle{(N_n)}}_{\scriptscriptstyle{\Lambda_{N_n}}}))$,
$n \in {\mathbb N}$, and ${\mathcal Q}$ the quadratic form corresponding to 
$({\mathcal E}^{\rm min}_\mu, D({\mathcal E}^{\rm min}_\mu))$.
In order to check (M1) we need to consider a closed
extension of $({\mathcal E}_\mu, {\mathcal F}C_b^{\infty}({\mathcal D}, \Gamma))$ 
on $L^2(\Gamma, \mu)$, which possibly is larger than
$({\mathcal E}^{\rm min}_\mu, D({\mathcal E}^{\rm min}_\mu))$. 
Let ${\mathcal VF}C_b^{\infty}({\mathcal D}, \Gamma)$ be the set of all maps defined 
as follows:
\begin{eqnarray*}
\Gamma \ni \gamma \mapsto \sum_{i=1}^N F_i(\gamma) v_i,
\end{eqnarray*} 
where $F_1, \ldots, F_N \in {\mathcal F}C_b^{\infty}({\mathcal D}, \Gamma)$ and
 $v_1, \ldots, v_N \in C^{\infty}_c({\mathbb R}^d, {\mathbb R}^d)$. 
For $V =  \sum_{i=1}^N F_i v_i \in {\mathcal VF}C_b^{\infty}({\mathcal D}, \Gamma)$
we define
\begin{eqnarray*}
{\rm div}_{\phi}^{\Gamma} V := \sum_{i=1}^N (\nabla^\Gamma_{v_i} F_i + B_{v_i}^{\phi} F_i v_i),
\end{eqnarray*} 
see \eqref{bofdiv} and Lemma \ref{pointwiselimit}. 
From the integration by parts formula \eqref{idibp} provided in 
Theorem \ref{thidibp} we can conclude that for all
$F \in {\mathcal F}C_b^{\infty}({\mathcal D}, \Gamma)$ and
$V \in {\mathcal VF}C_b^{\infty}({\mathcal D}, \Gamma)$    
\begin{eqnarray*}
\int_{\Gamma} \langle \nabla^\Gamma F, V \rangle_{T\Gamma} 
\,d\mu = - \int_{\Gamma} F
{\rm div}_\phi^\Gamma V \,d\mu.
\end{eqnarray*}
Let $(({\rm div}_\phi^\Gamma)^*, D(({\rm div}_\phi^\Gamma)^*))$ denote 
the adjoint of $({\rm div}_\phi^\Gamma, {\mathcal VF}C_b^{\infty}({\mathcal D}, \Gamma))$
as an operator from $L^2(\Gamma, T\Gamma, \mu)$ to $L^2(\Gamma, \mu)$. By
definition, $G \in L^2(\Gamma, \mu)$ belongs to $D(({\rm div}_\phi^\Gamma)^*)$ 
if and only if there exist a unique 
$({\rm div}_\phi^\Gamma)^* G \in L^2(\Gamma, T\Gamma, \mu)$ such that 
\begin{eqnarray*}
\int_{\Gamma} G \, {\rm div}_\phi^\Gamma V \,d\mu =
- \int_{\Gamma} \langle ({\rm div}_\phi^\Gamma)^* G, V \rangle_{T\Gamma} 
\,d\mu \quad \mbox{for all} 
\quad V \in {\mathcal VF}C_b^{\infty}({\mathcal D}, \Gamma). 
\end{eqnarray*}
We set $D({\mathcal E}^{\rm max}_{\mu}) :=
D(({\rm div}_\phi^\Gamma)^*)$, $d^\mu := ({\rm div}_\phi^\Gamma)^*$ and define
\begin{eqnarray*}
{\mathcal E}^{\rm max}_{\mu}(F,G) : =
\int_{\Gamma} \langle d^{\mu} F,
d^{\mu}G \rangle_{T\Gamma} \,d \mu, \quad \mbox{for all} \quad F,G 
\in D({\mathcal E}^{\rm max}_{\mu}).    
\end{eqnarray*}
From the integration by part formula \eqref{idibp} it follows that
\begin{eqnarray*}
{\mathcal F}C_b^{\infty}({\mathcal D}, \Gamma) \subset
D({\mathcal E}^{\rm max}_{\mu}) \quad \mbox{and} \quad d^\mu = \nabla^\Gamma 
\quad \mbox{on} \quad {\mathcal F}C_b^{\infty}({\mathcal D}, \Gamma).    
\end{eqnarray*}
Hence the densely defined, non-negative, closed, symmetric bilinear form
$({\mathcal E}^{\rm max}_{\mu},  D({\mathcal E}^{\rm max}_{\mu}))$    
extends $({\mathcal E}^{\rm min}_{\mu},  D({\mathcal E}^{\rm min}_{\mu}))$.
From general theory it is clear that $({\mathcal E}^{\rm max}_{\mu},  
D({\mathcal E}^{\rm max}_{\mu}))$ has an associated self-adjoint
generator $(H^{\rm max}_{\mu}, D(H^{\rm max}_{\mu}))$. However, it is not clear 
whether it is Markovian. In \cite{SchSu03} it is shown that for non-negative
interaction potentials $\phi$ the generator of the closure of 
${\mathcal E}^{\rm max}_{\mu}$ restricted to  
$D({\mathcal E}^{\rm max}_{\mu}) \cap L^\infty(\mu)$ is the maximum
Markovian self-adjoint extension of $( H_{\mu},
{\mathcal F}C_b^{\infty}({\mathcal D}, \Gamma))$.   

Condition (M1) we can only check, when ${\mathcal Q}$ is the
quadratic form corresponding to   
$({\mathcal E}^{\rm max}_{\mu}, D({\mathcal E}^{\rm max}_{\mu}))$.
Condition (M2) we can only check, when ${\mathcal Q}$ is the
quadratic form corresponding to   
$({\mathcal E}^{\rm min}_{\mu}, D({\mathcal E}^{\rm min}_{\mu}))$.
Hence we have to assume that $({\mathcal E}^{\rm min}_{\mu},  D({\mathcal E}^{\rm min}_{\mu}))
= ({\mathcal E}^{\rm max}_{\mu}, D({\mathcal E}^{\rm max}_{\mu}))$. In the case that
$(H^{\rm max}_{\mu}, D(H^{\rm max}_{\mu}))$ is Markovian, this is equivalent to
the so-called Markov uniqueness property, see e.g.~\cite{Eber99}. 
Obviously, the property 
$({\mathcal E}^{\rm min}_{\mu},  D({\mathcal E}^{\rm min}_{\mu}))$
$= ({\mathcal E}^{\rm max}_{\mu}, D({\mathcal E}^{\rm max}_{\mu}))$
is weaker than essential self-adjointness of $( H_{\mu},$
${\mathcal F}C_b^{\infty}({\mathcal D}, \Gamma))$. 

To verify (M1) we need strong convergence of the logarithmic derivatives.

\begin{proposition}\label{convld}
{\rm Let the conditions in Theorem \ref{tight} hold and 
$V \in {\mathcal VF}C_b^{\infty}({\mathcal D}, \Gamma)$.
Then ${\rm div}_{\phi}^{\Gamma} V$ considered as 
an element in $L^2(\mu^{\scriptscriptstyle{(N_n)}})$
converges strongly to ${\rm div}_{\phi}^{\Gamma} V $ 
considered as an element in $L^2(\mu)$ as $n \to \infty$
(recall that ${\rm div}_{\phi}^{\Gamma} V$ is a pointwise defined
function on $\Gamma$).}
\end{proposition}
{\bf Proof:} In the proof of Theorem \ref{thidibp} we have shown that
\begin{eqnarray}\label{weakly} 
\lim_{n \to \infty}
\int_{\Gamma_{\scriptscriptstyle{\Lambda_{N_n}}}^{\scriptscriptstyle{(N_n)}}} F
\, {\rm div}_\phi^\Gamma V \,d\mu^{\scriptscriptstyle{(N_n)}}
=  
\int_{\Gamma} F \, {\rm div}_\phi^\Gamma V \,d\mu
\end{eqnarray}
for all $F \in {\mathcal F}C_b^{\infty}({\mathcal D}, \Gamma)$. Hence, by
Lemma \ref{lemstrong} it remains to show that 
\begin{eqnarray}\label{eq844} 
\lim_{n \to \infty}
\int_{\Gamma_{\scriptscriptstyle{\Lambda_{N_n}}}^{\scriptscriptstyle{(N_n)}}} 
({\rm div}_\phi^\Gamma V)^2 \,d\mu^{\scriptscriptstyle{(N_n)}}
=  
\int_{\Gamma} ({\rm div}_\phi^\Gamma V)^2 \,d\mu.
\end{eqnarray}
Since for $V =  \sum_{i=1}^N F_i v_i \in {\mathcal VF}C_b^{\infty}({\mathcal D}, \Gamma)$
we have
\begin{eqnarray*}
{\rm div}_{\phi}^{\Gamma} V = \sum_{i=1}^N (\nabla^\Gamma_{v_i} F_i + 
\langle \mbox{div} \,v_i, \cdot \rangle + L_{v_i}^{\phi} F_i v_i),
\end{eqnarray*} 
\eqref{eq844} follows from 
\begin{align}\label{eq845} 
\lim_{n \to \infty}
\int_{\Gamma_{\scriptscriptstyle{\Lambda_{N_n}}}^{\scriptscriptstyle{(N_n)}}} 
|\langle \mbox{div} \,w, \cdot \rangle L_{v}^{\phi}| \,d\mu^{\scriptscriptstyle{(N_n)}}
& =  
\int_{\Gamma} |\langle \mbox{div} \,w, \cdot \rangle  L_{v}^{\phi}| \,d\mu \\
\mbox{and} \quad \lim_{n \to \infty}
\int_{\Gamma_{\scriptscriptstyle{\Lambda_{N_n}}}^{\scriptscriptstyle{(N_n)}}} 
(L_{v}^{\phi})^2 \,d\mu^{\scriptscriptstyle{(N_n)}}
& =  
\int_{\Gamma} (L_{v}^{\phi})^2 \,d\mu \label{eq846} 
\end{align}
for all $v, w \in C^{\infty}_c({\mathbb R}^d, {\mathbb R}^d)$. For all the other terms 
convergence can be shown as convergence of \eqref{weakly}. Using the
notation as in the proof of Lemma \ref{l2limit} we get:
\begin{multline*}
|{\mathbb E}_{\mu^{\scriptscriptstyle{(N_n)}}} [(L_v^{\phi})^2]
-{\mathbb E}_{\mu}[(L_v^{\phi})^2]| \le |{\mathbb
E}_{\mu^{\scriptscriptstyle{(N_n)}}}
[(L_v^{\phi})^2 - (L_{v, k}^{\phi})^2]| +
{\mathbb E}_{\mu^{\scriptscriptstyle{(N_n)}}}
[(L_{v,k}^{\phi})^2 - ((L_{v,k}^{\phi})^2)_r] \\
+ |{\mathbb E}_{\mu^{\scriptscriptstyle{(N_n)}}}
[((L_{v,k}^{\phi})^2)_r] - {\mathbb E}_{\mu}[((L_{v,k}^{\phi})^2)_r]|
+ {\mathbb E}_{\mu}[(L_{v,k}^{\phi})^2 - ((L_{v,k}^{\phi})^2)_r]
+ |{\mathbb E}_{\mu}[(L_{v,k}^{\phi})^2 - (L_{v}^{\phi})^2]| \\
\le \Big(\sup_{n \in {\mathbb N}}
\sqrt{{\mathbb E}_{{\mu}^{\scriptscriptstyle{(N_n)}}}[(L_{v}^{\phi})^2]}
+ \sup_{k \in {\mathbb N}}\sup_{n \in {\mathbb N}}
\sqrt{{\mathbb E}_{{\mu}^{\scriptscriptstyle{(N_n)}}}[(L_{v,k}^{\phi})^2]}
\Big)\sqrt{C_{16}(k)} \\ 
+ \Big(\sqrt{{\mathbb E}_{{\mu}}[(L_{v}^{\phi})^2]} + 
\sup_{k \in {\mathbb N}} \sqrt{{\mathbb E}_{{\mu}}[(L_{v,k}^{\phi})^2]}\Big)
\sqrt{{\mathbb E}_{\mu}[(L_v^{\phi} - L_{v, k}^{\phi})^2]} \\ +
|{\mathbb E}_{\mu^{\scriptscriptstyle{(N_n)}}}
[((L_{v,k}^{\phi})^2)_r] - {\mathbb E}_{\mu}[((L_{v,k}^{\phi})^2)_r]|
+ \frac{\sup_{k \in {\mathbb N}}
{\mathbb E}_{{\mu}}[|L_{v,k}^{\phi}|^3]}{r} \\
+ \frac{\sup_{k \in {\mathbb N}}\sup_{n \in {\mathbb N}}
{\mathbb E}_{{\mu}^{\scriptscriptstyle{(N_n)}}}[|L_{v,k}^{\phi}|^3]}{r}.
\end{multline*}
The constants $\sup_{k \in {\mathbb N}}{\mathbb
E}_{{\mu}}[|L_{v,k}^{\phi}|^3]$,
$\sup_{k \in {\mathbb N}}\sup_{n \in {\mathbb N}} {\mathbb
E}_{{\mu}^{\scriptscriptstyle{(N_n)}}}[|L_{v,k}^{\phi}|^3]$ are
finite due to Lemma \ref{l3norm}. 
Furthermore, the constants $\sup_{n \in {\mathbb N}}
{\mathbb E}_{{\mu}^{\scriptscriptstyle{(N_n)}}}[(L_{v}^{\phi})^2]$,
$\sup_{k \in {\mathbb N}}\sup_{n \in {\mathbb N}}
{\mathbb E}_{{\mu}^{\scriptscriptstyle{(N_n)}}}[(L_{v,k}^{\phi})^2]$
and 
$\sup_{k \in {\mathbb N}} {\mathbb E}_{{\mu}}[(L_{v,k}^{\phi})^2]$
are finite due to the estimates for $|J_2|, |J_3|, |J_4|$ provided in the
proof of Lemma \ref{l2limit}. Now 
\eqref{eq846} follows from Lemma \ref{l2limit} and the weak
convergence ${\mu}^{\scriptscriptstyle{(N_n)}} \to \mu$ as $n \to
\infty$. \eqref{eq845} can be shown analogously, since 
the constant $\sup_{n \in {\mathbb N}}
{\mathbb E}_{{\mu}^{\scriptscriptstyle{(N_n)}}}[|\langle \mbox{div} 
\,w, \cdot \rangle|^p]$ is finite for all $1 \le p < \infty$. 
For $p = 2$ this is shown in the proof of Theorem \ref{thidibp},
see \eqref{C17}. The proof easily generalizes to all $1 \le p < \infty$ due to the 
Ruelle bound \eqref{Ruellebo}. 
\hfill$\blacksquare$

\begin{theorem}\label{thmmoscoconv}
Let the assumptions in Theorem \ref{tight} hold and 
let $(\mu^{\scriptscriptstyle{(N_n)}
})_{n \in {\mathbb N}}$ be a subsequence
converging to an accumulation point $\mu$ of
$(\mu^{\scriptscriptstyle{(N)}})_{N \in {\mathbb N}}$.
Suppose that $({\mathcal E}^{\rm min}_{\mu},  D({\mathcal E}^{\rm min}_{\mu}))$
$= ({\mathcal E}^{\rm max}_{\mu}, D({\mathcal E}^{\rm max}_{\mu}))$.
Then the sequence of quadratic forms corresponding to
$({\mathcal E}^{\scriptscriptstyle{(N_n)}}_{\scriptscriptstyle{\Lambda_{N_n}}},$
$D({\mathcal E}^{\scriptscriptstyle{(N_n)}
}_{\scriptscriptstyle{\Lambda_{N_n}}}))_{n \in {\mathbb N}}$
Mosco converges to the quadratic form corresponding to
$({\mathcal E}^{\rm max}_{\mu},  D({\mathcal E}^{\rm max}_{\mu}))$.
\end{theorem}
{\bf Proof:} Since $({\mathcal E}^{\rm min}_{\mu}, D({\mathcal E}^{\rm min}_{\mu}))$
is the closure of $({\mathcal E}_\mu, {\mathcal F}C_b^{\infty}({\mathcal D}, \Gamma))$,
Proposition \ref{thmconvdf} together with \cite{Kol05}[Lem.~2.8] implies (M2) 

In order to check condition (M1), we consider 
a sequence $(F_n)_{n \in {\mathbb N}}$ with 
$F_n \in L^2(\mu^{\scriptscriptstyle{(N_n)}})$, $n \in {\mathbb N},$ which 
weakly converges to $F \in L^2(\mu)$. Furthermore, recall that in 
Proposition \ref{convld} we have shown strong convergence of 
${\rm div}_{\phi}^{\Gamma} V$ considered as 
an element in $L^2(\mu^{\scriptscriptstyle{(N_n)}})$
to ${\rm div}_{\phi}^{\Gamma} V $ 
considered as an element in $L^2(\mu)$ as $n \to \infty$ for 
all $V \in {\mathcal VF}C_b^{\infty}({\mathcal D}, \Gamma)$. 

First let us assume that
$F \in D({\mathcal E}^{\rm max}_\mu)$. Then
\begin{multline}\label{alex}
\Big( \int_{\Gamma} \langle d^\mu F, V \rangle_{T\Gamma} 
\,d\mu \Big)^2 = \Big(\int_{\Gamma} F
\, {\rm div}_\phi^\Gamma V \,d\mu\Big)^2 = \lim_{n \to \infty}
\Big(\int_{\Gamma_{\scriptscriptstyle{\Lambda_{N_n}}}^{\scriptscriptstyle{(N_n)}}} F_n
\, {\rm div}_\phi^\Gamma V \,d\mu^{\scriptscriptstyle{(N_n)}}\Big)^2 \\
\le \liminf_{n \to \infty} {\mathcal Q}_n(F_n) 
\int_{\Gamma} \langle V, V \rangle_{T\Gamma} \,d\mu \quad \mbox{for all} \quad
V \in {\mathcal VF}C_b^{\infty}({\mathcal D}, \Gamma),
\end{multline}
where ${\mathcal Q}_n$ is the quadratic form corresponding to 
$({\mathcal E}^{\scriptscriptstyle{(N_n)}}_{\scriptscriptstyle{\Lambda_{N_n}}},
D({\mathcal E}^{\scriptscriptstyle{(N_n)}}_{\scriptscriptstyle{\Lambda_{N_n}}}))$.
Since ${\mathcal VF}C_b^{\infty}({\mathcal D}, \Gamma)$ is dense in 
$L^2(\Gamma, T\Gamma, \mu)$, \eqref{alex} yields (M1) for $F \in D({\mathcal E}^{\rm max}_\mu)$.

For $F \notin  D({\mathcal E}^{\rm max}_\mu)$ we have $ {\mathcal Q}(F) = \infty$, where 
${\mathcal Q}$ is the quadratic form corresponding to 
$({\mathcal E}^{\rm max}_{\mu}, D({\mathcal E}^{\rm max}_{\mu}))$. We assume that 
$\liminf_{n \to \infty} {\mathcal Q}_n(F_n) = C_{19} < \infty$, then
\begin{eqnarray*}
\Big|\int_{\Gamma} F
\, {\rm div}_\phi^\Gamma V \,d\mu\Big| \le \sqrt{C_{19}}  \|V\|_{L^2(\Gamma, T\Gamma, \mu)} \quad 
\mbox{for all} \quad V \in {\mathcal VF}C_b^{\infty}({\mathcal D}, \Gamma).
\end{eqnarray*}
I.e., $F \in ({\rm div}_\phi^\Gamma)^* = D({\mathcal E}^{\rm max}_\mu)$. That is a contradiction! Hence
$\liminf_{n \to \infty} {\mathcal Q}_n(F_n) = \infty$.
\hfill$\blacksquare$

\begin{remark}
The essential ideas for proving Theorem \ref{thmmoscoconv} we got from 
the proof of \cite{Kol05}[Prop.~4.1].
\end{remark}

We denote the strongly continuous contraction semi-group associated with
$({\mathcal E}^{\rm max}_{\mu},$  $D({\mathcal E}^{\rm max}_{\mu}))$ by $(T_\mu(t))_{t \ge 0}$.

\begin{corollary}\label{convsemiresol}
Let the assumptions in Theorem \ref{tight} hold and 
let $(\mu^{\scriptscriptstyle{(N_n)}
})_{n \in {\mathbb N}}$ be a subsequence
converging to an accumulation point $\mu$ of
$(\mu^{\scriptscriptstyle{(N)}})_{N \in {\mathbb N}}$.
Suppose that $({\mathcal E}^{\rm min}_{\mu},  D({\mathcal E}^{\rm min}_{\mu}))$
$= ({\mathcal E}^{\rm max}_{\mu}, D({\mathcal E}^{\rm max}_{\mu}))$.
Then the sequence of semi-groups
$(T^{\scriptscriptstyle{(N_n)}}_{\scriptscriptstyle{\Lambda_{N_n}}}(t))$
strongly converges to $T_\mu(t)$ as $n \to \infty$ for all $t \ge 0$.
The same holds for the corresponding resolvents.  
\end{corollary}
{\bf Proof:} By Theorem \ref{thmmoscoconv} this follows directly from 
\cite[Theo.~2.4]{KuSh03}.
\hfill$\blacksquare$

\pagebreak

\begin{theorem}\label{uniqueness}
Let the assumptions in Theorem \ref{tight} hold and
let ${\bf P}$ be an accumulation point of $({\bf
P}^{\scriptscriptstyle{(N)}})_{N \in {\mathbb N}}$ with invariant
canonical Gibbs measure $\mu$. Suppose that 
 $({\mathcal E}^{\rm min}_{\mu},$  $D({\mathcal E}^{\rm min}_{\mu}))
= ({\mathcal E}^{\rm max}_{\mu}, D({\mathcal E}^{\rm max}_{\mu}))$.
Then ${\bf P}$ is the law of a 
Markov process with initial distribution $\mu$ 
and semi-group $(T_\mu(t))_{t \ge 0}$. In particular,
all accumulation points of 
$({\bf P}^{\scriptscriptstyle{(N)}})_{N \in {\mathbb N}}$ 
with the same invariant measure $\mu$ coincide.
\end{theorem}
{\bf Proof:} Let $({\bf P}^{\scriptscriptstyle{(N_n)}})_{n \in {\mathbb N}}$ 
be a subsequence such that $\lim_{n \to \infty } {\bf P}^{\scriptscriptstyle{(N_n)}} 
= {\bf P}$. This implies $\lim_{n \to \infty } {\mu}^{\scriptscriptstyle{(N_n)}} 
= {\mu}$. From Corollary \ref{convsemiresol} we now can conclude that 
$T^{\scriptscriptstyle{(N_n)}}_{\scriptscriptstyle{\Lambda_{N_n}}}(t)$
converges strongly to $T_\mu(t)$ as $n \to \infty$ for all $t \ge 0$.
Thus, finite dimensional distributions of ${\bf P}$ are 
given through $(T_\mu(t))_{t \ge 0}$. Since this holds for all 
accumulation points of $({\bf P}^{\scriptscriptstyle{(N)}})_{N 
\in {\mathbb N}}$ with invariant measure $\mu$, they all coincide.   
\hfill$\blacksquare$

\section{Application to the problem of equivalence of ensembles}\label{s5}

Grand canonical Gibbs measures correspond to an interaction
potential $\phi$, inverse temperature $\beta \ge 0$ and activity
function $z \ge 0$. An interesting question is the equivalence of
the grand canonical and canonical ensemble, i.e., the question
whether grand canonical and canonical Gibbs measures corresponding
to an interaction potential $\phi$ and inverse temperature $\beta$
coincide for a certain relation between their activity function
$z$ and particle density $\rho$, respectively, see
e.g.~\cite[Chap.~6]{Ge79}. Furthermore, it is of interest whether
one can approximate grand canonical Gibbs measures by finite
volume canonical Gibbs measures.

\begin{theorem}\label{coroequildht}
Assume that the conditions in Theorem \ref{tight} hold and
that we are in the low density, high temperature regime, i.e.,
\begin{eqnarray*}
\rho < \frac{1}{2\exp(2\beta K+1)J(\beta)}.
\end{eqnarray*}
Then the sequence of finite volume canonical Gibbs
measures with empty boundary condition
$(\mu^{\scriptscriptstyle{(N)}})_{N \in {\mathbb N}}$ converges to
a canonical Gibbs measure $\mu$ with constant density $\rho$
as $N \to \infty$. Furthermore, $\mu$ is a grand canonical
Gibbs measure corresponding to the activity
\begin{eqnarray}\label{activity}
z = \lim_{N \to \infty} N
\frac{Z^{\scriptscriptstyle{(N-1)}}_{\scriptscriptstyle{\Lambda_N}}}
{Z^{\scriptscriptstyle{(N)}}_{\scriptscriptstyle{\Lambda_N}}}.
\end{eqnarray}
\end{theorem}
{\bf Proof:} Let us fix an accumulation point 
$\mu$ of $(\mu^{\scriptscriptstyle{(N)}})_{N \in {\mathbb N}}$, i.e., 
there exists a subsequence $(\mu^{\scriptscriptstyle{(N_m)}}
)_{m \in {\mathbb N}}$ such that $\mu^{\scriptscriptstyle{(N_m)}} \to \mu$ 
weakly as $m \to \infty$. From Theorem \ref{thidibp} we know that
$\mu$ is a canonical Gibbs measure. Let $f \in C_c({\mathbb R}^d)$. Then
$\langle f, \cdot \rangle$ is a continuous function on $(\Gamma,
d_{\scriptscriptstyle{(\beta/3)}\Phi,h})$ and as in the proof of Theorem 
\ref{thidibp} we can show that 
\begin{eqnarray*}
{\mathbb E}_{\mu}[\langle f, \cdot \rangle] 
= \lim_{m \to \infty} {\mathbb E}_{\mu^{\scriptscriptstyle{(N_m)}}}
[\langle f, \cdot \rangle].
\end{eqnarray*}
Now \eqref{corre} yields 
\begin{eqnarray*}
{\mathbb E}_{\mu}[\langle f, \cdot \rangle] 
= \lim_{m \to \infty} 
\int_{\Lambda_{N_m}} f(x)\,
k_{\scriptscriptstyle{\Lambda_{N_m}}}^{\scriptscriptstyle{(1,N_m)}}
(x) \,dx_{\Lambda_{N_m}}.
\end{eqnarray*}
In \cite[Sect.~4]{BPK70} it is proved that in the low density, high temperature regime 
\begin{eqnarray*}
\lim_{m \to \infty} k_{\scriptscriptstyle{\Lambda_{N_m}}}^{\scriptscriptstyle{(1,N_m)}}
(x) = \rho  \quad \mbox{for all} \quad x \in {\mathbb R}^d. 
\end{eqnarray*}
Hence, using the Ruelle bound \eqref{Ruellebo} and Lebesgue's 
dominated convergence theorem, we obtain
\begin{eqnarray}\label{density}
{\mathbb E}_{\mu}[\langle f, \cdot \rangle] 
= \rho \int_{{\mathbb R}^d} f \,dx.
\end{eqnarray}
Thus, $\mu$ has constant density $\rho$. 

Moreover, in \cite[Sect.~4]{BPK70} it is proved that in the low density, high temperature regime
there exits a sequence $(k^{\scriptscriptstyle{(n)}})_{n \in {\mathbb N}}$ of functions
$k^{\scriptscriptstyle{(n)}}: {\mathbb R}^{n \cdot d} \to {\mathbb R}$ such that
\begin{eqnarray*}
\lim_{m \to \infty} k_{\scriptscriptstyle{\Lambda_{N_m}}}^{\scriptscriptstyle{(n,N_m)}}
(x_1, \ldots, x_n) = k^{\scriptscriptstyle{(n)}}(x_1, \ldots, x_n) \,\, \mbox{for all} 
\,\, (x_1, \ldots, x_n) \in {\mathbb R}^{n\cdot d} \,\, \mbox{and all} \,\,
n \in {\mathbb N}. 
\end{eqnarray*}  
Using the Ruelle bound, see Theorem \ref{ruelle}, and analogous arguments as in the 
derivation of \eqref{density}, we can identify the sequence 
$(k^{\scriptscriptstyle{(n)}})_{n \in {\mathbb N}}$ as the correlation functions of $\mu$.
Since this is true for all accumulation points of $(\mu^{\scriptscriptstyle{(N)}})_{N \in {\mathbb N}}$,
all accumulation points coincide and  
$(\mu^{\scriptscriptstyle{(N)}})_{N \in {\mathbb N}}$ converges to 
the canonical Gibbs measure $\mu$ as $N \to \infty$.   

Finally, \cite[Theo.~I]{BPK70} together with \cite[Theo.~IV]{BPK70}
implies that the sequence of correlation functions
$(k^{\scriptscriptstyle{(n)}})_{n \in {\mathbb N}}$
fulfills the Kirkwood--Salsburg equations for
$z$ as given in \eqref{activity}. Thus, $\mu$ is a grand
canonical Gibbs measure corresponding to 
$\phi$, $\beta$ and $z$. 
\hfill$\blacksquare$

\begin{remark}\label{microensemble}
{\rm A related result has been proved in \cite{Geo95}. There the
author derived an approximation of grand canonical Gibbs measures
by finite volume micro canonical Gibbs measures with periodic
boundary condition. It seems to be quite feasible to adapt the
proof to the canonical case. However, since that proof heavily
relays on the choice of a periodic boundary condition, it still
would not cover our case (empty boundary condition).}
\end{remark}

\begin{corollary}
Assume that the conditions in Theorem \ref{tight} hold and
that we are in the low density, high temperature regime.
Let $\mu = \lim_{N \to \infty} \mu^{\scriptscriptstyle{(N)}}$
and suppose that $({\mathcal E}^{\rm min}_{\mu},$  
$D({\mathcal E}^{\rm min}_{\mu})) = ({\mathcal E}^{\rm max}_{\mu}, 
D({\mathcal E}^{\rm max}_{\mu}))$. Then the sequence $({\bf
P}^{\scriptscriptstyle{(N)}})_{N \in {\mathbb N}}$ converges in
law to a Markov process ${\bf P}$ with initial distribution $\mu$ 
and semi-group 
$(T_\mu(t))_{t \ge 0}$. Furthermore, ${\bf P}$ solves the martingale 
problem for $(-H_{\mu}, {\mathcal F}C_b^{\infty}({\mathcal D}, \Gamma))$ 
with initial distribution $\mu$. 
\end{corollary}
{\bf Proof:} This is an immediate consequence of Theorem \ref{uniqueness}
together with Theorem \ref{coroequildht} and Theorem \ref{martingaleproblem}.
\hfill$\blacksquare$

\addcontentsline{toc}{section}{References}


\end{document}